\documentclass{amsart}
\usepackage{fullpage}
\usepackage{tikz} % package p r i n c i p a l TikZ
\usetikzlibrary{arrows} % l i b r a i r i e o p t i onne l l e PGF
\usepackage{xcolor}
\usepackage[english]{babel}
\usepackage{epsfig}
\usepackage{latexsym}
\usepackage{amsmath}
\usepackage{amssymb}

\newtheorem{coro}{Corollary}
\newtheorem{thm}{Theorem}
\newtheorem{lem}{Lemma}
\newtheorem{prop}{Proposition}

\newcommand{\bbR}{\mathbb{R}}
\newcommand{\bbN}{\mathbb{N}}

\newcommand{\norm}{\|}
\newcommand{\bvert}{\boldsymbol{\vert}}
\newcommand{\lb}{\boldsymbol{[\![}}
\newcommand{\rb}{\boldsymbol{]\!]}}
\newcommand{\un}{\boldsymbol{u}}
\newcommand{\hun}{\boldsymbol{hu}}
\newcommand{\p}{\boldsymbol{p}}
\newcommand{\q}{\boldsymbol{q}}

\begin{document}
\title{Max-Plus Convexity in Riesz spaces}
\author{Charles Horvath}
\address{Universit\'e de 
Perpignan D\'epartement de Math\'ematiques, Laboratoire LAMPS 
52 av. Paul Alduy 66860  PERPIGNAN cedex 9
}
\thanks{horvath@univ-perp.fr}
\begin{abstract}
We study max-plus convexity in an Archimedean Riesz space $E$ with an order unit $\un$; the definition of max-plus convex sets is algebraic and we do not assume that $E$ has an {\it a priori} given topological structure. To the given unit $\un$ one can associate two equivalent  norms $\norm\cdot\norm_{\un}$ and $\norm\cdot\norm_{\hun}$ on $E$; the distance 
${\sf D}_{\hun}$ on $E$ associated to $\norm\cdot\norm_{\hun}$ is  a geodesic distance for which  max-plus convex sets in $E$ are  geodesically closed sets. Under suitable assumptions, we establish max-plus versions of some   fixed points and continuous selection theorems that are well known for linear convex sets and we    show that hyperspaces of compact max-plus convex sets are Absolute Retracts. \end{abstract}

\maketitle

 \section{Introduction} 
 To keep the size of this paper reasonable, the definition of max-plus convexity in arbitrary Riesz spaces with respect to a given unit $\un$ is given in Section \ref{defmax+conv} with barely no justification as to why one should be interested in max-plus convexity. If the Riesz space in question is $\bbR^n$ and $\un = (1, \cdots, 1)$ then the max-plus convex sets with respect to $\un$ are the usual max-plus convex subsets of $\bbR^n$. The usual finite dimensional max-plus convexity lives in $\big(\bbR \cup\{-\infty\}\big)^n$, the extension to arbitrary Archimedean Riesz spaces with a unit that is presented here is therefore not a full generalization of the finite dimensional framework which does not mean that such a thing could not be done. The reader looking for motivations and applications is referred to \cite{work},  \cite{tropalfgeo},  \cite{litmas}, \cite{propgeoams}. 
 Section \ref{prelim} is mainly about  basic  concepts  and a few exemples. 
 Section \ref{tropnorm} is about two norms,  denoted here by $\norm \cdot \norm_{\un}$ and   $\norm \cdot \norm_{\hun}$, which can be associated to a given unit $\un$ in an Archimedean Riesz space $E$. If $E = \bbR^n$ and 
 $\un = (1, \cdots, 1)$ then   $\norm (x_1, \cdots, x_n) \norm_{\un} = \max\vert x_i\vert$ and 
  $\norm (x_1, \cdots, x_n) \norm_{\hun} = \max x_i^+ + \max x_i^-$,  the so called Hilbert affine norm. In Section \ref{defmax+conv} one will find the definition of max-plus convex sets in an Archimedean Riesz space, with respect to a given unit, and some of their basic properties the most important one being the Kakutani Property, also known in the standard linear framework as the Algebraic Hahn-Banach Property. Section \ref{geodesics} shows that ${\sf D}_{\hun}$, the metric associated to $\norm \cdot \norm_{\hun}$,  is a geodesic distance on the Riesz space $E$ with respect to which the geodesically closed sets are precisely the max-plus convex sets, with respect to the given unit $\un$. That max-plus convexity in $\bbR^n$ should be a geodesic structure with respect to an appropriate metric is a not so rescent idea; it had been discussed years ago with Walter Briec from Universit\'e de Perpignan and it is Stefan Gaubert, in a discussion with the author at \'Ecole Polytechnique, who hinted at the fact that the Hilbert affine metric should be the appropriate metric; shortly thereafter, in a private communication \cite{gauperso} mailed to the author, Stefan Gaubert proved that this is indeed the case.  The proof given here differs somewhat from Gaubert's straightforward coordinatewise proof in $\bbR^n$ but would have been impossible without that proof. Section \ref{topmax} deals with the basic topological properties of max-plus convex sets and hyperspaces of compact max-plus sets;   Section \ref{topmax+conv} is about infinite dimensional max-plus versions of some standard results: Ky Fan Best Approximation and conequently Schauder's Fixed Point Theorem, Kakutani's Fixed Point Theorem for upper semicontinuous maps, Michael's Selection Theorem, Dugundji's Extension Theorem, and consequently the fact that max-plus convex sets are Absolute Retracts.   Max-plus convexity in $\bbR^S$ and hyperspaces of compact max-plus convex subsets of $\bbR^S$ have been studied by  L. Bazylevych, D. Repovs and  M. Zarichnyi in \cite{bazyrepza}. Hyperspaces of max-plus compact convex sets in $\mathcal{C}(X)$, where $X$ is a compact metrizable topological space, have been studied by  L. Bazylevych and  M. Zarichnyi in \cite{bazyzar}.  In Section \ref{topmax} one can find a few remarks on hyperspaces of compact max-plus convex sets in a Riesz space .

\section{Preliminaries and a few examples}\label{prelim}
We will denote by $\bbR_+$ the set of positive real numbers and by $\bbR_{++}$ the set of strictly positive real numbers. 

\medskip
A {\bf Riesz space}, or a {\bf vector lattice},  is a real vector space $E$ endowed with a partial order $\leqslant$ that is compatible with the linear structure, that is 
\begin{equation}\label{defrieszspa}
\forall x, y, z \in E\, \, \forall t\in\bbR_+\quad x\leqslant y \Rightarrow (tx + z)\leqslant (ty+z)\\
\end{equation}  and such that all pairs $\{x, y\}$ of elements of $E$ have a least upper bound for which we will use the standard notation $x\vee y$.

The positive cone is $E_+ = \{x\in E : 0\leqslant x\}$ which has the following properties; 
\begin{equation}\label{prope+}
\begin{cases}
(1)& \forall x, y\in E \quad x\leqslant y \Leftrightarrow (y-x)\in E_+\\
(2)& \forall\, \, x, y\in E_+ \, \forall t\in\bbR_+\, \, (tx +y)\in E_+
\end{cases}
\end{equation}

All  pairs $\{x, y\}$ of elements of $E$ have a greatest lower bound, for which we will use the standard notation  
$x\wedge y$; one easily sees that $x\wedge y =  -\big((-x)\vee (-y)\big)$. 
 
 \medskip A Riesz space $E$ is {\bf  Archimedean}     if,  whenever $x$ and $y$ are  two elements of $E$ such that, for all $n\in \bbN$, $ny\leqslant x$, one has $y \leqslant 0$. 

\medskip  A Riesz space $E$ is {\bf Dedekind complete} (respectively, {\bf Dedekind $\boldsymbol{\sigma}$-complete}) if every non-empty (respectively, countable) subset $S$ of $E$ which has an upper bound has a least upper bound. \\
  Since $S$ has an upper bound if and only if $-S$ has a lower bound one can replace  in the definition of Dedekind complete (resp. $\sigma$-completness) ``upper bound'' by ``lower bound''. A Dedekind $\sigma$-complete Riesz space is Archimedean.

  \medskip\noindent Every Archimedean Riesz space has a Dedekind completion, more precisely: there exists a Dedekind complete Riesz space $\hat{E}$ containing $E$ as a vector sublattice such that
  \begin{equation}\label{archdensdede}
  \forall \hat{x}\in\hat{E}\quad\hat{x} = \sup\{x\in E : x\leqslant\hat{x}\} =  \inf\{x\in E : \hat{x}\leqslant x\}
  \end{equation}

 A {\bf strong order unit} of a Riesz space $E$ is an element $\un\in E_+$ such that
\begin{equation}\label{defstrgu} 
\forall x\in E_+\, \,  \exists n\in\bbN \text{ such that }  x\leqslant n \un 
\end{equation}

Since strong units are the only kind of units we will consider  we will drop the adjective ``strong''.

  \bigskip A {\bf Riesz norm} on a Riesz space $E$ is a norm  such that, 
\begin{equation}\label{defRnorm}
\forall x, y\in E\quad  
 \bvert x \bvert \leqslant \bvert y \bvert \Rightarrow \norm x \norm \leq \norm y \norm
 \end{equation} 
 
 For a Riesz norm one has, for all $x\in E$,     $\norm\bvert x \bvert \norm = \norm x \norm$. \footnote{$x^+ = 0\vee x$,  
 $x^- = -(0\wedge x) = 0\vee(-x)$ and $ \bvert x \bvert = x\vee (-x) = x^+ + x^- = x^+\vee x^-$. Also,  $\bvert tx\bvert = \vert t \vert\, \bvert x \bvert$,  $\bvert x + y \bvert \leqslant \bvert x \bvert + \bvert y \bvert$, 
$\bvert x \bvert \leqslant y$ if and only if $-y \leqslant x \leqslant y$, 
$(x+ y)^+ \leqslant x^+ +y^+$,  $(x+ y)^- \leqslant x^-+ y^-$, 
$(x\vee y)^+ = x^+ \vee y^+$, $(x\vee y)^- = x^- \wedge y^-$.}
 
 \medskip\noindent A Riesz space equipped with a Riesz norm is a {\bf normed lattice}. A normed lattice is Archimedean and the lattice operations are uniformly continuous.
 
 \medskip Let $S$ be a subset of a Riesz space  $E$ for which there exist $x_1, x_2\in E$ such that, for all $x\in S$, $x_1\leqslant x \leqslant x_2$ ($S$ is an order bounded set); if $\norm \cdot \norm$ is a Riesz norm on $E$ then, 
 $x\in S$, $\norm x_1 \norm\leqslant \norm x \norm \leqslant \norm x_2 \norm$. That is,  in a normed lattice an order bounded set is norm-bounded.  

 \medskip\noindent Any two complete lattice norms on a given Riesz space $E$ are equivalent,  page 352 or  \cite{frem} Proposition 25 A.
  
   \bigskip An {\sf M}{\bf -norm} on a Riesz space $E$ is a Riesz norm such that, 
   \begin{equation}\label{defMnorm}
   \forall x, y\in E_+\quad\norm x\vee y\norm = \max\{\norm x\norm, \norm y \norm\}
   \end{equation}

\medskip\noindent If $\norm \cdot \norm$ is an {\sf M}-norm on  $E$ and if $0\leqslant x\leqslant y$ then $\norm x \norm \leq \norm y \norm$ ( from $x\vee y = y$, 
 $\norm  y\norm = \max\{\norm x\norm, \norm y \norm\}$). If  $\norm \cdot \norm$ is an {\sf M}-norm and a Riesz norm on  $E$ then, for all $x, y\in E$, $\norm x\vee y\norm \leq \max\{\norm x\norm, \norm y \norm\}$.

 \medskip\noindent An {\sf AM}{\bf -space} is a Riesz space equipped with a complete norm which is an {\sf M}-norm and a  Riesz norm. 
 
 \medskip\noindent An  {\sf AM}{\bf -space with a unit} (resp.  an {\sf M}{\bf -space with a unit}) is an  {\sf AM}-space $E$ 
 (resp.  an {\sf M}-space ) with a unit $\un$ such that $\norm \un \norm = 1$ in which case the unit ball  is 
 $\{x\in E : -\un \leqslant x \leqslant\un\}$.

 \bigskip There is a standard way to associate to each given unit $\un$ on an Archimedean Riesz space $E$ an {\sf M}-norm  
 $\norm \cdot \norm_{\un}$  on $E$, namely  : 
 \begin{equation}\label{defstdnrm}
 \norm x \norm_{\un} = \inf\{t\in\bbR_\star :\,  \bvert x \bvert \leqslant t\un\} 
 \end{equation} 
 
  \medskip That norm $\norm \cdot \norm_{\un}$ is  a Riesz norm is evident; that it is also an {\sf M}-norm is well known, it can also be seen from  Lemma \ref{pqmaxmim} below. 
  
   \medskip If $\un_1$ and $\un_2$ are two units on an Archimedean Riesz space $E$ then, the norms 
 $\norm \cdot \norm_{\un_{1}}$ and $\norm \cdot \norm_{\un_{2}}$ are equivalent since there exists two natural numbers 
 $n_1$ and $n_2$ such that $\un_1\leqslant n_1\un_2$ and $\un_2\leqslant n_2\un_1$.

  \medskip In an Archimedean Riesz space $E$, possibly without a unit, take an arbitrary element  $\un\in E_+\!\setminus\!\{0\}$ and let \\
 $E_{\un} = \{x\in E : \exists n\in \bbN \, \, \bvert x \bvert \leqslant n\un\}$ (the principle ideal spanned by $\un$). Then  $\big(E_{\un},\norm\cdot\norm_{\un}\big)$ is an {\bf M}-space with unit. If $\un$ is a unit of $E$ then $E_{\un}$ is $E$ itself. Furthermore,  if $(E, \norm \cdot \norm )$ is a complete normed lattice or   if $E$ is Dedekind $\sigma$-complete   then, for all $\un\in E$,  $\big(E_{\un},\norm\cdot\norm_{\un}\big)$ is an {\bf AM}-space with unit; details can be found  \cite{frem}, Lemma 25I and Lemma 25J. 

\medskip If the Riesz space $E$ is equipped with a complete Riesz norm $\norm \cdot \norm$ and if  $\un$ is a unit in $E$ then $\norm \cdot \norm$ and $\norm \cdot \norm_{\un}$ are complete lattice norms on $E$; they  are therefore equivalent.

 \bigskip If $\Omega$ is a compact topological space then the space $\mathcal{C}(\Omega)$ of continuous real valued functions on $\Omega$ is a Riesz space, the fucntion $\boldsymbol{1}$ identically equal to $1$ is a unit and, for all 
 $x\in \mathcal{C}(\Omega)$, $\norm x \norm_{\boldsymbol{1}} = \sup_{\omega\in\Omega}\vert x(\omega)\vert = \norm x \norm_{\infty}$; $\big( \mathcal{C}(\Omega), \norm \cdot \norm_{\infty}\big)$ is of course an  {\bf AM}-space with unit.

\bigskip  The {\bf Bohnenblust-Kakutani-Krein Representation Theorem}   says that ${\sf AM}$-spaces with a unit are isomorphic, as normed Riesz spaces\footnote{Riesz spaces $E_1$ and $E_2$ are isomorphic if there exists a linear isomorphism which is also a lattice isomorphism.  },  to spaces of continuous functions on a compact topological space $\Omega$; see \cite{luxza}, Chapter 17, \&121. 
  
 \bigskip\noindent 
 {\bf (BKK RepresentationTheorem)} {\it An {\sf AM}-space with  unit is Riesz isomorphic and norm isomorphic to a space $\mathcal{C}(\Omega)$ equipped with the $\sup$-norm,  for an appropriate compact Hausdorff space $\Omega$.}
 
\bigskip Now, consider an Archimedean Riesz space $E$ with unit $\un$ and let $\hat{E}$ be its Dedekind completion; by the preceding discussion, 
$\big(\hat{E}_{\un}, \norm \cdot \norm_{\un}\big)$ is an {\sf AM}-space with unit; by the BBK-Representation Theorem there exits  a compact topological space $\Omega$ and a Riesz space isomorphism  $\Phi : \hat{E}_{\un}\to \mathcal{C}(\Omega)$ such that, for all $\hat{x}\in \hat{E}_{\un}$, $\norm \hat{x} \norm_{\un} = \norm \Phi(\hat{x}) \norm_{\infty}$ and 
$\Phi(\un) = \boldsymbol{1}$;   since $E$ is itself a Riesz subspace of $\hat{E}_{\un}$ and $\un\in E$; $E$ can  be identified  with the Riesz subspace $\Phi(E) = \tilde{E}$ of $\mathcal{C}(\Omega)$. In conclusion, 

\medskip\noindent{\it An Archimedean Riesz $E$ space with unit $\un$ can be identified with a Riesz subspace of a space 
$\mathcal{C}(\Omega)$ where $\Omega$ is a compact Hausdorff space that is, there is a Riesz subspace ${\tilde E}$ of  $\mathcal{C}(\Omega)$ such that $\boldsymbol{1}\in {\tilde E}$ and there exists a Riesz space isomorphism 
$\Phi : E\to {\tilde E}$ such that, for all $x\in E$, $\norm x \norm_{\un} = \norm \Phi(x) \norm_{\infty}$.}

\medskip The closure  $\bar{E}$ of the Riesz space $\tilde{E}$ in $\mathcal{C}(\Omega)$   is an 
{\sf AM}-space, and since  $\boldsymbol{1}\in {\tilde E}\subseteq \bar{E}$ it is an {\sf AM}-space with unit, since  $\boldsymbol{1}$ is a unit in $\mathcal{C}(\Omega)$,    which shows that an Archimedean Riesz space possessing a  unit can be embedded as a dense Riesz subspace of an {\sf AM}-space with unit and therefore as dense Riesz subspace of some  $\mathcal{C}(\Omega^\prime)$ containing $\boldsymbol{1}$, where $\Omega^\prime$ is some compact topological space. If the Riesz space $\tilde{E}$  separates the points of $\Omega$ then $\bar{E} = \mathcal{C}(\Omega)$, by the  Stone-Weierstrass Theorem.

\bigskip We conclude this section with  a few simple examples of the previous constructions.

\medskip For all sets $S$, the space $\mathcal{F}(S) = \bbR^S$ of arbitrary real valued functions on $S$  with pointwise operations, is an Archimedean Riesz space, without a unit if $S$ is not a finite set; $\mathcal{F}_\star(S)$, the space of bounded real valued functions on $S$, is an Archimedean Riesz space with a unit: $\un = (1, 1, 1, \ldots)$.\\
Let $\un\in \mathcal{F}(S)$ be a positive function which is not identically $0$. For all $x\in \mathcal{F}(S)$ let $Z(x) = \{\omega\in S : x(\omega) = 0\}$. Let $E = \mathcal{F}(S)$; then  $x\in E_{\un}$ if $Z(x)\subset Z(\un)$ and $\displaystyle{\sup_{\omega\not\in Z(\un)} \frac{\vert x(\omega)\vert }{\un(\omega)} < \infty}$ in which case $\norm x \norm_{\un}$ is this supremum.  
If  $\un : S \to\bbR$ is identically equal to one the $ E_{\un} = \mathcal{F}_\star(S)$.

\medskip  If $\Omega$ is a compact topological space then ${\mathcal C}(\Omega)$ is Dedekind complete (resp.  Dedekind $\sigma$-complete) if and only if the closure of every open set (resp. of every $F_\sigma$ open set) is open.  If   $\un\in\mathcal{C}(\Omega)$ is a unit, that is a strictly positive function, then 
 \begin{equation*}
 \norm x \norm_{\un} =  \displaystyle{\max_{\omega\in \Omega}\frac{\vert x(\omega) \vert}{\un(\omega)}}.
 \end{equation*}

\medskip The classical sequence spaces $l_p$ are Riesz subspaces of $\bbR^\bbN$; the sequence whose terms are all equal to $1$ is a unit of $l_\infty$.

\medskip  Let  $\big(X, \mathcal{B}, \mu\big)$ be a measured space where the measure $\mu$ is  finite; ${\mathcal M}\big(X, \mathcal{B}, \mu\big)$, the space of almost 
$\mu$- everywhere finite  real valued functions, with the usual identification of almost $\mu$-everywhere equal functions,   is a  Dedekind complete Riesz space.  Also, still under the hypothesis that $\mu$ is  a finite measure,   the spaces $L_p\big(X, \mathcal{B}, \mu\big)$, $1\leq p\leq \infty$ are Dedekind complete. Details can be found in \cite{luxza1} page $126-127$.

 \medskip\noindent Given a measured space $\left(\Omega, \mathcal{B}, \mu\right)$, $L_\infty\left(\Omega, \mathcal{B}, \mu\right)$ is also, with for its usual norm, an {\sf AM}-space:  the constant map $\un(\omega) = 1$ is a unit.

\medskip 
Let $\Omega$ be a non empty set, $\mathcal{F}$ a field of subsets of $\Omega$ and 
${\sf ba}(\Omega, \mathcal{F})$ the family of bounded charges\footnote{$\emptyset\in \mathcal{F}$  on $\mathcal{F}$ and $\Omega\in \mathcal{F}$ ; if $A, B\in \mathcal{F}$ then $A\cup B\in\mathcal{F}$ and $A\setminus\! B\in\mathcal{F}$. An element of ${\sf ba}(\Omega, \mathcal{F})$ is a  map $\mu : \mathcal{F}\to\bbR$ which is additive ($\mu(A\cup B) = \mu(A) + \mu(B)$ if $A\cap B = \emptyset$) and such that $\mu(\emptyset) = 0$ and $\sup\vert\mu(A)\vert < \infty$.} on $\mathcal{F}$; 
${\sf ba}(\Omega, \mathcal{F})$ is a real vector space and the relation $\mu \leqslant \nu$ if, for all $A\in\mathcal{F}$, $\mu(A) \leq \nu(A)$, is a partial order on ${\sf ba}(\Omega, \mathcal{F})$ compatible with the vector space structure;  endowed with that partial order ${\sf ba}(\Omega, \mathcal{F})$ is a Riesz space, the maximum of two elements $\mu$ and $\nu$ of 
${\sf ba}(\Omega, \mathcal{F})$ is given by, for all $A\in\mathcal{F}$,  
\begin{equation*}
(\mu \vee\nu) (A) = \sup\big\{\mu(B) + \nu(A\setminus\!B) : B\subseteq A \text{ and } B\in\mathcal{F} \big\}.
\end{equation*}
Furthermore, ${\sf ba}(\Omega, \mathcal{F})$ is Dedekind complete and $\mu\mapsto \bvert\mu\bvert(\Omega)$ is a complete Riesz norm on ${\sf ba}(\Omega, \mathcal{F})$ but not an {\sf M}-norm\footnote{By definition of the absolute value of an element of a Riesz space, $\bvert\mu\bvert = \mu^+ + \mu^-$;  an explicit formula for the absolute value of 
a bounded charge $\mu$ is $\bvert\mu(A)\bvert = \sup_{\mathcal{R}}\big\{\sum_{F\in\mathcal{R}}\vert\mu(F)\vert : \mathcal{R}\text{ is a finite partition of } A \text{ by elements of } \mathcal{F}\big\}$.}; by the discussion above, if $\un\in{\sf ba}(\Omega, \mathcal{F})$ is a positive charge then  ${\sf ba}(\Omega, \mathcal{F})_{\un}$ equipped with the norm 
$\norm \cdot \norm_{\un}$ is an {\sf AM}-space. 

\medskip A charge $\mu\in {\sf ba}(\Omega, \mathcal{F})$ is a measure if,  for all countable family $\{F_n : n\in\bbN\}$ of pairwise disjoint elements of $ \mathcal{F}$ whose union belongs  to $ \mathcal{F}$, one has $\mu(\cup_{n}F_n) = \sum_{n}\mu(F_n)$; the set of elements of ${\sf ba}(\Omega, \mathcal{F})$ that are measures is a Dedekind-complete Riesz subspace of ${\sf ba}(\Omega, \mathcal{F})$ which is also a Banach sublattice. Details can be found in \cite{rao}.

   \bigskip We assume throughout that $\boldsymbol{E}$ {\bf is an Archimedean Riesz space} and that $\un$ {\bf is a unit of} $\boldsymbol{E}$.   

\section{Max-plus norms in Archimedean spaces}\label{tropnorm} 
 For all $x\in E$ let 
 \begin{equation}\label{defqupu}
 \begin{cases}
 \p_{\un}(x) = \inf\{t\in\bbR : x\leqslant   t\un\}\\ 
\q_{\un}(x) = \sup\{t\in\bbR : t\un  \leqslant x\}\\
 \end{cases}
 \end{equation}
 
\medskip\noindent One has 
 \begin{equation}\label{p-q-}
 \p_{\un}(x) = -\q_{\un}(-x)
 \end{equation}
 and 
 \begin{equation}\label{monot}
 \forall x, y\in E\quad x\leqslant y\Rightarrow\begin{cases}
\p_{\un}(x)\leq \p_{\un}(y)
\\ \q_{\un}(x)\leq \q_{\un}(y)\\
  \end{cases}
 \end{equation}
 
 \medskip Since $\un$ is a unit, $\{t\in\bbR : x\leqslant   t\un\}\neq\emptyset$ from which $\p_{\un}(x) < \infty$. If $\p_{\un}(x)$ were $-\infty$ then we would have, for all 
 $t\in\bbR$, $t\un \leqslant (-x)$ and, since $E$ is Archimedean, we would have $\un\leqslant 0$; which is not the case. In conclusion, $\p_{\un}(x)$ is a real number and, by (\ref{p-q-}), so is $\q_{\un}(x)$.

  \begin{lem}\label{qxp}
  \begin{equation*}
  \forall x\in E\quad\q_{\un}(x)\un\leqslant x \leqslant\p_{\un}(x)\un
  \end{equation*}
  
  from which $x = 0$ if and only if $\q_{\un}(x) = \p_{\un}(x) = 0$. 
  \end{lem}
  
  \proof If we show that $ x \leqslant\p_{\un}(x)\un$ then the other the inequality will follow from (\ref{qxp}).\\
  If $\p_{\un}(x) < t$ then $x \leqslant t\un$ and therefore, for all $n\in\bbN_\star$, $x\leqslant (\p_{\un}(x) + \frac{1}{n})\un$ that is $n(x-\p_{\un}(x)\un)\leqslant\un$ and since $E$ is Archimedean, $x - \p_{\un}(x)\un\leqslant 0$. $\Box$
  
  \bigskip Call a $\un$-box in $E$ any order interval of the form $[s\un, t\un] = \{y\in E : s\un\leqslant y \leqslant t\un\}$, where $s\leq t$ are real numbers; then, by definition, $[\q_{\un}(x), \p_{\un}(x)]$ is the smallest $\un$-box containing $x$.
  
 \begin{lem}\label{pqmaxmim}
 \begin{equation*}
 \forall x, y\in E \begin{cases}
(1)\quad\p_{\un}(x\vee y) =\max\{\p_{\un}(x), \p_{\un}(y)\}\\
 (2)\quad\q_{\un}(x\wedge y) = \min\{\q_{\un}(x), \q(y)\}
 \end{cases}
 \end{equation*}
 From which we have
 
 \begin{equation*}
 \forall x\in E\quad \p_{\un}(x^+) = \max\{0, \p_{\un}(x)\} \text{ and } \p_{\un}(x^-) = \max\{0, \p_{\un}(-x) \}
  \end{equation*}
  
  and also $x =0$ if and only if $\p_{\un}(x^+) = \p_{\un}(x^-) = 0$. 
 \end{lem} 
 
 \proof From $x\leqslant\p_{\un}(x)\un$ and  $y\leqslant\p_{\un}(y)\un$ we have $x\vee y \leqslant\max\{\p_{\un}(x), \p_{\un}(y)\}\un$ and therefore 
 $\p_{\un}(x\vee y)\leq \max\{\p_{\un}(x), \p_{\un}(y)\}$.\\
 From $x\leqslant x\vee y$ and $y\leqslant x\vee y$ we have $\p_{\un}(x)\leqslant \p_{\un}(x\vee y)$ and $\p_{\un}(y)\leqslant \p_{\un}(x\vee y)$ from which we have  $ \max\{\p_{\un}(x), \p_{\un}(y)\}\leq \p_{\un}(x\vee y)$. 
 
 \medskip From (\ref{p-q-}) and $x\wedge y = - \big((-x)\vee (-y)\big)$ we have $\q_{\un}(x\wedge y) = \min\{\q_{\un}(x), \q_{\un}(y)\}$.
 
 \medskip From $x^+ = 0\vee x$ we have $\p_{\un}(x^+) = \max\{\p_{\un}(0), \p_{\un}(x)\}$ and, from 
 $x^- = 0\vee (-x)$ we have \\
 $\p_{\un}(x^-) = \p_{\un}(0\vee (-x)) =  \max\{0, \p_{\un}(-x)\}$. 
 
 \medskip If $\p_{\un}(x^+) = 0 = \p_{\un}(x^-)$  then $\p_{\un}(x)\leq 0$ and $\p_{\un}(-x)\leq 0$; from $x\leqslant \p_{\un}(x)\un$, 
 $-x\leqslant\p_{\un}(-x)\un$ and $\un\in E_+$ we have $x\in E_+$ and $-x\in E_+$ and therefore $x = 0$.
 $\Box$
 
\begin{lem}\label{homgsub}
\begin{equation*}
\forall x, y\in E\text{ and } \forall s\in\bbR_+\,\, \begin{cases}
(1)\, 
\p_{\un}(x + y)\leq \p_{\un}(x) + \p_{\un}(y) \text{ and } \p_{\un}(sx) = s\p_{\un}(x)\\
(2)\,\q_{\un}(x) + \q_{\un}(y)\leq \q_{\un}(x + y) \text{ and } \q_{\un}(sx) = s\q_{\un}(x)
\end{cases}
\end{equation*}
and, if $s < 0$ then $\p_{\un}(sx) = s\q_{\un}(x)$. 
\end{lem}

\proof From $x\leqslant\p_{\un}(x)\un$ and  $y\leqslant\p_{\un}(y)\un$ we have $x + y \leqslant(\p_{\un}(x) + \p_{\un}(y))\un$ from which the first part of $(1)$ follows.  For the second part, there is nothing to prove if $s = 0$; if $s > 0$ then 
$sx\leqslant(s\p_{\un}(x))\un$ from which we have $\p_{\un}(sx)\leq s\p_{\un}(x)$. \\
From $sx\leq\p_{\un}(sx)\un$ we have $x\leq(s^{-1}\p_{\un}(sx))\un$ and therefore, $\p_{\un}(x)\leq s^{-1}\p_{\un}(sx)$ that is 
 $s\p_{\un}(x)\leq \p_{\un}(sx)$.
 
 \medskip The second part is a consequence of $\q_{\un}(x) = -\p_{\un}(-x)$. 
 
 \medskip If $s < 0$ then $\p_{\un}(sx) = \p_{\un}\big((-s)(-x)\big) = \vert s \vert \p_{\un}(-x) = (-\vert s \vert)(-\p_{\un}(-x)) = s\q_{\un}(x)$. $\Box$

\begin{thm}\label{hunrienrm} The map from $E$ to $\bbR$ given by  
\begin{equation*}
x\mapsto \p_{\un}(x^+) + \p_{\un}(x^-) = \norm x \norm_{\hun}
\end{equation*} 
is a norm   on $E$. Furthermore, for all $x, y\in E$, $\norm x\vee y\norm_{\hun}\leq\max\{\norm x\norm_{\hun}, \norm y\norm_{\hun}\}$ with equality if both $x$ and $y$ are in $E_+$ and
 \begin{equation}\label{nrm++}
 \forall x\in E\quad \norm x \norm_{\hun} = \norm x^+ \norm_{\hun} + \norm x^- \norm_{\hun} =  \norm x^+ \norm_{\un} + \norm x^- \norm_{\un} 
 \end{equation}
\end{thm}

\proof From $x^+\in E_+$ and $x^-\in E_+$  we have $\p_{\un}(x^+)\geq 0$ and $\p_{\un}(x^-)\geq 0$ and therefore 
$\norm x \norm_{\hun}\geq 0$. If $\norm x \norm_{\hun} = 0$ then $\p_{\un}(x^+) = \p_{\un}(x^-) =  0$ and by Lemma 
\ref{pqmaxmim}, $x = 0$. Clearly,  $\norm 0 \norm_{\hun} = 0$. 

\medskip If $s > 0$ the $(sx)^+ = sx^+$ and $(sx)^- = sx^-$; from  Lemma \ref{homgsub} we obtain 
$\norm sx \norm_{\hun} = s\norm x \norm_{\hun}$. 

\medskip If $s < 0$ then $sx = \vert s\vert(-x)$ from which $\norm sx \norm_{\hun} = \vert s\vert \, \norm -x \norm_{\hun}$ and, $\norm -x \norm_{\hun} = \p_{\un}((-x)^+) + \p_{\un}((-x)^-) = \p_{\un}(x^-) + \p_{\un}(x^+) = \norm x \norm_{\hun}$.  

\medskip $\norm x + y\norm_{\hun}\leq \norm x \norm_{\hun} + \norm y\norm_{\hun}$ follows From 
$(x+y)^+ \leqslant x^+ + y^+$, $(x+y)^- \leqslant x^- + y^-$ and Lemma \ref{homgsub} we have 
$\norm x + y\norm_{\hun} = \p_{\un}\big(\big(x+y)^+\big) +  \p_{\un}\big((x+y)^-\big) \leqslant \p_{\un}(x^+) + \p_{\un}(y^+) + \p_{\un}(x^-) + \p_{\un}(y^-) = 
\norm x \norm_{\hun} + \norm  y\norm_{\hun}$.

\medskip Without loss of generality let us assume that $\p_{\un}(x^+)\leq \p_{\un}(y^+)$. From $(x\vee y)^+ = x^+\vee y^+$
 and $(x\vee y)^- = x^-\wedge y^-$ and from $(1)$ of Lemma \ref{pqmaxmim} we have 
 $\norm x\vee y\norm_{\hun} \leq \max\{\p_{\un}(x^+), \p_{\un}(y^+)\} + \min\{\p_{\un}(x^-) + \p_{\un}(y^-)\}\leq 
 \p_{\un}(y^+) + \p_{\un}(y^-)$. Finally, if $x$ and $y$ are both in $E_+$ then  $x\vee y \in E_+$. $\Box$

 \bigskip We will call $\norm \cdot \norm_{\hun}$ the {\bf max-plus norm} on $E$ associated to the unit $\boldsymbol{u}$. 
 
 \medskip From Theorem \ref{hunrienrm},  the max-plus norm $\norm \cdot \norm_{\hun}$ is an {\sf M}-norm, but it is not a Riesz norm as can be seen by taking $E = \bbR^2$, $\un = (1, 1)$ and $x = (-1, 1)$ for which we have $\norm x \norm_{\hun} = 2$  and $\norm\bvert x \bvert \norm_{\hun} = 1$.\\ More generally, we always have 
 $\norm \bvert x \bvert \norm_{\hun} \leq \norm x \norm_{\hun}$ since  $\norm \bvert x \bvert \norm_{\hun} = 
 \norm x^+ + x^-\norm_{\hun}\leq \norm x^+\norm_{\hun} + \norm x^-\norm_{\hun} = \p_{\un}(x^+) + \p_{\un}(x^-)$.

 \begin{coro}\label{defhilbdist}
 The map $(x, y)\mapsto {\sf D}_{\hun}(x, y) =   \max\{0, \p_{\un}(x-y)\} - \min\{0, \q_{\un}(x-y)\}$\\
 $\phantom{i}$\\
 $\phantom{AAAAAAAAAAAAAAAAAAAAAAAAA} = \max\{0, \p_{\un}(x - y)\} + \max\{0, \p_{\un}(y-x)\}$\\ is a  metric on $E$.
 \end{coro}
 
 \bigskip We will call ${\sf D}_{\hun}$ the {\bf max-plus distance} on $E$ associated to the unit $\boldsymbol{u}$.

 \bigskip  For $E = \bbR^n$, $\un = (1, \cdots, 1)$ and $x = (x_1, \cdots, x_n)$ we have\footnote{$[n] = \{1, \cdots, n\}$} 
 
 \medskip
 \centerline{$\begin{cases}
 \p_{\un}(x) =  \inf\{t\in\bbR : \forall i\in [n] \, \, x_i\leq t \} = \max\{x_i : i\in [n]\}\\
 \q_{\un}(x) =  \sup\{t\in\bbR : \forall i\in[n] \, \, t\leq x_i\} = \min\{x_i : i\in [n]\}
 \end{cases}$}
 
 \medskip\noindent
 From $\p_{\un}(-x) = -\min\{x_i : i\in [n]\} = \max\{-x_i : i\in [n]\}$ we finally have 
 \begin{center}$\begin{cases}\max\{0, \p_{\un}(x)\} = \max\{x^+_i : i\in [n]\} = \p_{\un}(x^+)\\   
 \max\{0, \p_{\un}(-x)\} = \max\{x^-_i : i\in [n]\} = \p_{\un}(x^-)
 \end{cases}$\end{center}
  which gives 
 \begin{equation}\label{classichilbnorm}
 \norm x \norm_{\hun} =   \max\{x^+_i : i\in [n]\} + \max\{x^-_i : i\in [n]\} = \p_{\un}(x^+) + \p_{\un}(x^-)
\end{equation}

Put simply, if $x\not\in\bbR^n_+$ or if $-x\not\in\bbR^n_+$ then $ \norm x \norm_{\hun}$ is the difference between the largest and the smallest coordinate of $x$; if $x\in\bbR^n_+$, $ \norm x \norm_{\hun}$ is the largest coordinate of $x$ and, 
if $-x\in\bbR^n_+$, $ \norm x \norm_{\hun}$ is the absolute value of the smallest coordinate of $x$ that is, in this last two case, $ \norm x \norm_{\hun} = \norm x \norm_\infty = \max_{i\in [n]}\vert x_i \vert$. The norm $ \norm \cdot \norm_{\hun}$ on $\bbR^n$ associated to the unit $\un = (1, \cdots, 1)$ is the {\bf Hilbert affine norm} of S. Gaubert. 
   
  \bigskip Let $\boldsymbol{U}\big(E\big)$ be the set of units   of the Riesz space $E$. The set of units of  
$\bbR^n$ is $\bbR^n_{++} = \{x\in\bbR^n_+ : 0 < \min_{i\in[n]}x_i\}$. 

\begin{lem}\label{equivhunun}
For all $\boldsymbol{u}_1, \boldsymbol{u}_2\in \boldsymbol{U}\big(E\big)$ the norms  
$\norm \cdot \norm_{\boldsymbol{hu}_i}$ and $\norm \cdot \norm_{\boldsymbol{u}_j}$, $i, j\in\{1, 2\}$ are equivalent.

\end{lem}
\proof First, take $\boldsymbol{u}_1 = \boldsymbol{u}_2 = \boldsymbol{u}$. 
From the definitions,   $\norm x \norm_{\un} = \p\big(\bvert x \bvert\big)$ and therefore  
 $\norm x \norm_{\un} = \norm x \norm_{\hun}$ for all $x\in E_+$ which gives $\norm x \norm_{\hun} =
  \norm x^+ \norm_{\un} + \norm x^- \norm_{\un} \leq 2 \norm x \norm_{\un} $.

  \medskip\noindent From $\norm x\norm_{\un} = \norm \bvert x\bvert \norm_{\un}$ we get  
  $\norm x\norm_{\un} =  
\norm x^+  + x^-\norm_{\un}\leq \norm x^+\norm_{\un} + \norm x^-\norm_{\un} = \p(x^+) + \p(x^-) = \norm x\norm_{\hun}$. 

\medskip To complete the proof, notice that, for  $\boldsymbol{u}_1, \boldsymbol{u}_2\in \boldsymbol{U}\big(E\big)$, the norms $\norm \cdot \norm_{\boldsymbol{hu}_1}$ and $\norm \cdot \norm_{\boldsymbol{hu}_2}$ are equivalent since  
  $n_1 \boldsymbol{u}_2 \leqslant \boldsymbol{u}_1\leqslant n_2\boldsymbol{u}_1$ for some non zero whole numbers $n_1$ and $n_2$. 
$\Box$

\bigskip From Lemma \ref{equivhunun} one has that $\norm \cdot \norm_{\boldsymbol{hu}}$ is complete for  a given $\boldsymbol{u}\in  \boldsymbol{U}\big(E\big)$ if and only if $\norm \cdot \norm_{\boldsymbol{hu}}$ is complete for  all $\boldsymbol{u}\in  \boldsymbol{U}\big(E\big)$ if and only if $\norm \cdot \norm_{\boldsymbol{u}}$ is complete for  a given $\boldsymbol{u}\in  \boldsymbol{U}\big(E\big)$ if and only if $\norm \cdot \norm_{\boldsymbol{u}}$ is complete for  all $\boldsymbol{u}\in  \boldsymbol{U}\big(E\big)$.

\begin{lem}\label{contlatop}
For all $\boldsymbol{u}\in  \boldsymbol{U}\big(E\big)$ the lattice operations $\vee$, $\wedge$ are uniformly continuous with respect to the metric ${\sf D}_{\hun}$. Furthermore, if $K\subset E$ is ${\sf D}_{\hun}$-compact then there exists $t_1, t_2\in\bbR$ such that, for all $x\in K$, $t_1\un \leqslant  x \leqslant t_2\un$.
\end{lem} 
 
\proof The norms  $\norm \cdot \norm_{\un}$ and  $\norm \cdot \norm_{\hun}$  are equivalent and $\norm \cdot \norm_{\un}$ is a lattice norm, which implies uniform continuity of the lattice operations with respect to $\norm \cdot \norm_{\un}$.
$\Box$

\section{Max-plus convexity in Archimedean Riesz spaces with a  unit}\label{defmax+conv}
Given a Riesz space $E$  and a unit   $\un$   of $E$ one can, as in the now standard finite dimensional case 
$E = \bbR^n$,  
 introduce on $E$ two operations, {\bf max-plus addition} $\oplus$ and {\bf max-plus multiplication} $\odot$ {\bf by real numbers}: 

\begin{equation}\forall (x, y, t)\in E\times E\times\bbR\quad
\begin{cases}
x\oplus y\,  \text{ for }\,  x\vee y\,  \text{ and}\\
t\odot x\,  \text{ for }\,   x + t\un.
\end{cases}
\end{equation}

\medskip\noindent It would have been more appropriate to write $t\odot_{\un} x$  for $x + t\un$ since this ``tropical multiplication'' depends on the chosen unit $\un$. 

\medskip Furthermore, writing  ${\sf 1}$ for the real number $0$, $(\bbR, \odot, {\sf 1})$ is a group (the additive group of real numbers written multiplicatively) and  
$(t, x)\to t\odot x$ is an action of the group  $(\bbR, \odot, {\sf 1})$ on $E$ with the following properties: 

\begin{equation}\label{basicmp}
\begin{cases}
(1)&t\odot(x\oplus y) = (t\odot x)\oplus (t\odot y)\\
(2)& (t_1\oplus t_2)\odot x = (t_1\odot x) \oplus (t_2\odot x)\\
(3)&t_2\odot(t_1\odot x) = (t_1\odot t_2)\odot x\\
(4)&{\sf 1}\odot x = x
\end{cases}
\end{equation}

The notation makes everything look very familiar; the peculiarity here is that the ``sum'' is idempotent: $x\oplus x = x$. 

\medskip In the finite dimensional case one can enlarge the set of scalars  to $\bbR\cup\{-\infty\}$; tropical addition and multiplication are extended to $\bbR\cup\{-\infty\}$ in the obvious way: $(-\infty)\oplus x  = x\oplus(-\infty)= x$ and $(-\infty)\odot x = -\infty$. Writing ${\sf 0}$ for $-\infty$ this becomes ${\sf 0}\oplus x = x\oplus {\sf 0} = x$ and 
${\sf 0}\odot x = {\sf 0}$.

\medskip $E = \big(\bbR^n\cup\{-\infty\}\big)^n$ with pointwise $\oplus$ and $\odot$ operations, is an indempotent semimodule over the idempotent semi-field $\bbR\cup\{-\infty\}$; this is the standard  max-plus semi-module $\bbR_{max+}^n$.

 \bigskip Such a coordinatewise  procedure that turns an arbitrary Riesz space with unit into an idempotent semi-module over the semi-field 
 $\big({\bbR}\cup\{-\infty\}, \oplus, \odot, {\sf 0}, {\sf 1})$ is not  readily available. By  adding a single element to $E$ which becomes by decree the smallest element, the Riesz space $E$ is embedded in a semi-module over $\bbR\cup\{-\infty\}$ of which it is a max-plus convex subset(the definition of max-plus convexity is given below).   
 
 \medskip To an arbitrary Riesz space $E$  add a smallest  element, let us call it $\boldsymbol{\bot}$, and one can extend scalar multiplication to $\odot : \bbR\cup\{-\infty\}\times\big(E\cup\{\boldsymbol{\bot}\}\big) \to \big(E\cup\{\boldsymbol{\bot}\}\big)$ and  the tropical sum  to $ \big(E\cup\{\boldsymbol{\bot}\}\big)\times \big(E\cup\{\boldsymbol{\bot}\}\big)$ in such a way that $E\cup\{\boldsymbol{\bot}\}$ becomes an idempotent semimodule over the totally ordered semi-field $\bbR\cup\{-\infty\}$:

 \begin{equation}\label{extendop}
 \begin{cases}
 (1)\quad \forall x\in E\cup \{\boldsymbol{\bot}\}& \boldsymbol{\bot}\oplus x = x\oplus \boldsymbol{\bot} = x\\
 (2)\quad\forall x\in E\cup \{\boldsymbol{\bot}\}& -\infty\odot x = \boldsymbol{\bot}\\
 (3)\quad \forall t\in \bbR\cup\{-\infty\}& t\odot \boldsymbol{\bot} = \boldsymbol{\bot}
  \end{cases}
 \end{equation}
 
The {\bf max-plus convex hull} of a     non empty subset $S\subset E\cup \{\boldsymbol{\bot}\}$ is the set of elements  of $E\cup \{\boldsymbol{\bot}\}$ which can be written as $( x_1 + t_1\un)\vee\cdots\vee (x_m + t_m\un)$  with 
 $\{x_1, \cdots, x_m\}\subset S$, $\{t_1, \cdots, t_m\}\subset\bbR\cup\{-\infty\}$ and  $\max\{t_1, \cdots, t_m\} = 0$\footnote{In tropical notation this becomes $(t_1\odot x_1)\oplus\cdots\oplus (t_m\odot x_m)$ with $t_1\oplus \cdots \oplus t_m = {\sf 1}$
 (recall here ${\sf 1}$ is the usual $0\in\bbR$ and  that $\odot$ depends on a fixed unit $\un$ : $t\odot x = x + t\un$) exactly as if it were a usual affine combination.}.   We will use the notation $\lb S \rb_{\un}$ for the   max-plus convex hull of $S$ with respect to the unit $\un$; we set 
 $\lb \emptyset \rb_{\un} = \emptyset$. Whenever a single fixed unit $\un$ is under consideration we drop the index $\un$.
 
 \medskip\noindent
 If $S\subset E$ then $\lb S \rb_{\un}$ is a subset of $E$; it is the set of of elements of $E$ which can be written as\\ $( x_1 + t_1\un)\vee\cdots\vee (x_m + t_m\un)$ with 
 $\{x_1, \cdots, x_m\}\subset S$, $\{t_1, \cdots, t_m\}\subset\bbR$ and  $\max\{t_1, \cdots, t_m\} = 0$ since, if $t_i = -\infty$ then $ x_i + t_i\un = \boldsymbol{\bot}$ and, since one the coefficients is $0$, let us say $t_j = 0$, we have 
 $( x_i + t_i\un)\vee( x_j+ t_j\un) = \boldsymbol{\bot}\vee x_j = x_j$.

\medskip A subset $C$ of $E\cup\{-\infty\}$ is said to be {\bf max-plus convex} (with respect to $\un$) if $C = \lb C \rb_{\un}$. \\ 
 If $S = \{x_1, \cdots, x_m\}$ is a finite set we will write $\lb x_1, \cdots, x_m\rb$ for $\lb S \rb$. The max-plus convex hull of two points $x_1$ and $x_2$ is $\lb x_1, x_2\rb = \{x_1\vee(x_2 + t\un) : t\leq 0\} \cup \{(x_1 + t\un)\vee x_2 : t\leq 0\}$; if $x_1\neq x_2$, the set $\lb x_1, x_2\rb$ will be called a {\bf max-plus segment}. Since $\un$ is a unit, there exists $s \geq 0$ such that $(x_2-x_1)\leqslant s\un$ that is $x_2 + (-s)\un\leqslant x_1$ which shows that $x_1\in\lb x_1, x_2\rb$, and similarly for $x_2$.
 
 \medskip Given a non empty set $S$ we will denote by $\langle S \rangle$ the family of non empty finite subsets of $S$. A set of the form $\lb S\rb$ with $S$ finite is a {\bf max-plus polytope} (with respect to $\un$); this definition makes the empty set into a polytope.

 \medskip Max-plus convex subsets of $E\cup\{\boldsymbol{\bot}\}$ will be rarely referred to. Unless otherwise specified, ``max-plus convex set'' will mean ``max-plus convex set of $E$''.   What matters here, is that $E$ is a max-plus convex subset of the $(\bbR\cup\{-\infty\})$-idempotent semimodule $E\cup\{\boldsymbol{\bot}\}$ and therefore, the max-plus convex subsets of $E$ are exactly the max-plus convex subsets of $E\cup\{\boldsymbol{\bot}\}$ that are contained in $E$.
 
\medskip The proof of Lemma \ref{mainconvlem1} below is left to the reader;  that of Lemma \ref{recur} can here be done by hand or one can go to  Lemma 2.1.5 in \cite{moietamm}.
   
 \begin{lem}\label{mainconvlem1}
 Given a Riesz space $E$ a unit $\un$ of $E$ the following properties hold:\\
 $(1)$ $S \subset \lb S \rb$.\\
 $(2)$ If $S_1\subset S_2$ then $\lb S_1 \rb\subset \lb S_2 \rb$. \\
 $(3)$ $\displaystyle{\lb S \rb = \bigcup_{A\in\langle S \rangle}\lb A \rb}$ where $\langle S \rangle$ denotes the set of non empty finite subsets of $S$. \\
 $(4)$ $ \lb\, \lb S \rb\, \rb = \lb S \rb$.
 \end{lem}
 
  \begin{lem}\label{recur} 
For all finite subset $S$ of $E$ and for all $x\in E$, 
\begin{equation*}
\displaystyle{\lb S \cup \{x\} \rb = \bigcup_{y\in\lb S \rb}\lb x, y\rb}
\end{equation*}
\end{lem}

\begin{lem}\label{eqsup}$\phantom{i}$\\ 
(1) For all $x_1, x_2\in E$, $\lb x_1, x_2\rb = \lb x_1, x_1\vee x_2\rb \cup \lb x_2, x_1\vee x_2\rb$ and $\lb x_1, x_1\vee x_2\rb \cap \lb x_2, x_1\vee x_2\rb = \{x_1\vee x_2\}$. \\
(2)  If $w\in\lb x_1, x_2\rb$ then either $w\vee x_2 = x_1\vee x_2$ or $w\vee x_2 = w$.\\
(3) If $x_1$ and $x_2$ are comparable then $\lb x_1, x_2\rb$ is, with respect to the partial order of $E$, a totally ordered subset.
 \end{lem}
  
   \proof (1) If $z = x_1\vee(x_2 + s\un)$ with $s\leq 0$ then, from $(x_1 + s\un)\vee(x_2 + s\un) = (x_1\vee x_2) + s\un$ and $x_1 + s\un \leqslant x_1$ we have $x_1\vee \big( (x_1\vee x_2) + s\un\big) = x_1\vee (x_1 + s\un)\vee(x_2 + s\un) = x_1\vee(x_2 + s\un) = z$ which shows that $z\in \lb x_1, x_1\vee x_2\rb$; similarly, if 
   $z = (x_1 + s\un)\vee x_2$ with $s\leq 0$ then $z\in \lb x_2, x_1\vee x_2\rb$.\\
   From $\{x_i, x_1\vee x_2\}\subset \lb x_1, x_2\rb$ and from Lemma \ref{mainconvlem1} we have 
   $\lb x_i, x_1\vee x_2\rb\subset \lb x_1, x_2\rb$. \\
   Take $z\in \lb x_1, x_1\vee x_2\rb \cap \lb x_2, x_1\vee x_2\rb$ ; if $z = tx_i\vee (x_1\vee x_2)$ with $t\leq 0$ then $z = x_1\vee x_2$. If $z = x_1\vee s(x_1\vee x_2) = x_2\vee t(x_1\vee x_2)$ with $s, t\leq 0$ then 
   $z = z\vee z = x_1\vee _2\vee s(x_1\vee x_2) \vee s(x_1\vee x_2) = x_1\vee x_2$. This completes the proof of (1).\\
    (2)   If $w \in \lb x_1, x_1\vee x_2\rb$ then $x_1\leqslant w\leqslant x_1\vee x_2$ and therefore $x_1\vee x_2 \leqslant w\vee x_2 \leqslant  x_1\vee x_2\vee x_2 =  x_1\vee x_2$. \\
    If $w \in\lb x_2, x_1\vee x_2\rb$ then $x_2\leqslant w\leqslant x_1\vee x_2$ and therefore 
   $w\vee x_2 = w$. \\
   (3) Assume that $x_1\leqslant x_2$. An arbitrary element of $\lb x_1, x_2\rb$ is either of the form $(x_1 + s\un)\vee x_2$ or $x_1\vee (x_2 + s\un)$ with $s\leq 0$; since $(x_1 + s\un)\vee x_2 = x_2$ and 
   $x_1\leqslant x_1\vee (x_2 + s\un) \leq x_1\vee x_2 = x_2$, $x_1$ is the smallest element of $\lb x_1, x_2\rb$ and $x_2$ is its largest element. Take $z$ and $z^\prime$ in $\lb x_1, x_2\rb$; we can assume that $z = x_1\vee (x_2 + s\un)$ and $z^\prime = x_1\vee (x_2 + s^\prime\un)$ with $s\leq 0$ and $s^\prime\leq 0$. We have either 
   $s\leq s^\prime$ of $s^\prime\leq s$ from which it follows that either $z\leqslant z^\prime$ or 
   $z^\prime\leqslant z$.
    $\Box$

\begin{lem}\label{mainconvlem2}$\phantom{i}$\\
$(1)$ A singleton is max-plus convex. \\
$(2)$ An arbitrary intersection of max-plus convex sets is max-plus convex.\\
$(3)$ $\lb S \rb$ is the smallest max-plus convex set containing $S$. \\
$(4)$ $C$ is max-plus convex if and only if, for all $x_1, x_2\in C$, $\lb x_1, x_2\rb \subset C$.\\
$(5)$ The following statements are equivalent: 
 
 \begin{equation}\label{easyC}
 \begin{cases}
 (a)\,  C \text{ is max-plus convex}.\\
 (b)\,  \forall t \leq 0\quad\forall x, y\in C\quad x\vee (y + t\un)\in C\\
 (c)\, \big[\, \forall x_1, x_2\in C\quad x_1\vee x_2\in C\big] \text{ and }  \big[\, \forall x_1, x_2\in C \text{ such that } 
x_1\leqslant x_2\quad  \lb x_1, x_2\rb \subset C\,\big].
 \end{cases}
 \end{equation}
 \end{lem}

\proof Only $(4)$ needs to be checked, $(1), (2)$ and $(3)$ are direct consequences of the definitions. Let us see that $C$ is max-plus convex if for all $x_1, x_2\in C$, $\lb x_1, x_2\rb \subset C$; the reverse implication is a consequence of $(2)$ and $(4)$ of Lemma \ref{mainconvlem1}. By $(3)$ of Lemma  \ref{mainconvlem1} we have to show that, for all $A\in \langle C \rangle$, $\lb A\rb \subset C$ which can be proved using Lemma \ref{recur} and an obvious induction on the cardinality of the set $A$. \\
The equivalence of $(4)$ and $(5)$ follows from $\lb x_1, x_2\rb = \lb x_1, x_1\vee x_2\rb \cup  \lb x_2, x_1\vee x_2\rb$. $\Box$

 \bigskip A {\bf max-plus half-space} of a convex subset $C$ of $E$ is a max-plus convex set $D\subset C$ such that $C\setminus\! D$ is also max-plus convex.  A max-plus convex subset $C$ of $E$ has the {\bf Kakutani-Property} if, for arbitrary disjoint max-plus convex subsets $C_1$ and $C_2$ there exists a half-space $D$ of $C$ such that $C_1\subset D$ and $C_2\cap D = \emptyset$.
 
 \begin{lem}\label{kakuprop}
 Max-plus convex subsets of $E$ have the Kakutani-Property.
 \end{lem}

\proof From 3.0.12 in \cite{moietamm}. $\Box$

\bigskip We conclude this section by showing that the max-plus metric is additive on max-plus segments,\\ Proposition \ref{geod1} below. 

\begin{lem}\label{split}
For $x_1, x_2\in E$ 
\begin{equation*}
{\sf D}_{\hun}(x_1, x_2) = {\sf D}_{\hun}(x_1, x_1\vee x_2) + {\sf D}_{\hun}(x_2, x_1\vee x_2).
\end{equation*}
\end{lem}

 \proof If $x_1$ and $x_2$ are comparable then either $x_1 = x_1\vee x_2$ or $x_2 = x_1\vee x_2$ in which case there is nothing to prove. If  $\p_{\un}(x_1-x_2)\leq 0$ then $x_1\leqslant x_2$ and $x_1$ and $x_2$ are comparable, and similarly if $\p_{\un}(x_2-x_1)\leq 0$.  We can therefore assume that 
 $0 \leq \p_{\un}(x_1-x_2)$ and  $0 \leq \p_{\un}(x_2-x_1)$.

\medskip We now have  ${\sf D}_{\hun}(x_1, x_2) =  \max\{0, \p_{\un}(x_1-x_2)\} + \max\{0, \p_{\un}(x_2-x_1)\} = \p_{\un}(x_1-x_2) + \p_{\un}(x_2-x_1)$ 
 and  we have to see that $\p_{\un}(x_1-x_2) + \p_{\un}(x_2-x_1)\geq \norm x_1 -  x_2\vee x_1\norm_{\un} + \norm x_2 -  x_2\vee x_1\norm_{\un}$.

\medskip From $x_1 \leqslant \p_{\un}(x_1-x_2)\un + x_2$ and $\p_{\un}(x_1-x_2)\un\in E_+$ we have 
$x_1\vee x_2 \leqslant \p_{\un}(x_1-x_2)\un + x_2$ and, since  $\norm \cdot \norm_{\un}$ is a Riesz norm and 
 $\norm \un \norm_{\un} = 1$, $\norm x_1\vee x_2 - x_2\norm_{\un}\leq  \p_{\un}(x_1-x_2)$. Similarly, 
 $\norm x_1\vee x_2 - x_1\norm_{\un}\leq  \p_{\un}(x_2-x_1)$. 
 $\Box$

  \begin{lem}\label{wxypu}
  If $x_1\leq x_2$ then, for all $z\in\lb x_1, x_2\rb$, ${\sf D}_{\hun}(z, x_2) =  \p_{\un}(x_2- z) = \norm x_2- z\norm_{\un}$,\\ $z = x_1\vee\big(x_2 - \p_{\un}(x_2- z)\un\big)$ and 
  $\q_{\un}(z - x_2) = \sup\{s\leq 0 : z = x_1\vee(x_2 + s\un)\}$.  
  \end{lem}
  
  \proof From $z\in\lb x_1, x_2\rb$ we have either $z = (x_1 + s\un)\vee x_2$ or $z = x_1\vee(x_2 + s\un)$ with $s\leq 0$. In the first case, from $x_1\leq x_2$, we have $z = x_2$ in which case  $\p_{\un}(x_2- z) = 
  \p_{\un}(0) = 0 = \q_{\un}(0)$ and $x_1\vee\big(x_2 - \p_{\un}(x_2- z)\un\big) = x_1\vee x_2 = x_2 = z$.
  
  \medskip If $z = x_1\vee(x_2 + s\un)$ with $s\leq 0$ then $x_1\vee z = z$ and $x_2\vee z = x_2$, 
     which shows that 
  $x_1\leqslant z \leqslant x_2$. 
  
  \medskip\noindent From $0\leqslant x_2 - z$ we have $\p_{\un}(x_2- z) = \inf\{t\geq 0 :   x_2 - z \leqslant t\un\} = {\sf D}_{\hun}(z, x_2) = \norm x_2- z\norm_{\un}$.
  
  \medskip\noindent From $x_2 + s\un\leqslant x_1\vee(x_2 + s\un) = x_1\vee z = z$ we have 
  $\p_{\un}(x_2- z) \leq (-s)$ from which we obtain\\  
  $x_2 + s\un \leqslant x_2-\p_{\un}(x_2- z)\un$ and $z = x_1\vee(x_2 + s\un)\leqslant x_1\vee \big(x_2 - \p_{\un}(x_2- z)\un\big)$.
  
  \medskip\noindent From Lemma \ref{qxp} we have $x_2- z \leqslant \p_{\un}(x_2- z)\un$ that is, 
  $x_2 - \p_{\un}(x_2- z)\un \leqslant z$ and, since $x_1\leqslant z$, \\
  $x_1\vee \big(x_2 - \p_{\un}(x_2- z)\un\big) \leqslant z$.
  
  \medskip We have shown that $z = x_1\vee \big(x_2 - \p_{\un}(x_2- z)\un\big)$ and that 
  $s \leq -\p_{\un}(x_2- z) = \q_{\un}(z - x_2)$.  $\Box$

     \begin{lem}\label{splitbis} If $x_1$ and $x_2$ are two comparable elements of $E$ then, for all 
 $z\in \lb x_1, x_2\rb$, 
   \begin{equation}\label{additivity1}{\sf D}_{\hun}(x_1, x_2) = {\sf D}_{\hun}(x_1, z) + {\sf D}_{\hun}(z, x_2).
 \end{equation}
 \end{lem}
 
  \proof   From Lemma \ref{wxypu} we can write    $z = x_1\vee (x_2 + s\un)$ with $s = -{\sf D}_{\hun}(z, x_2)$. \\ From $x_1\leqslant z$,  ${\sf D}_{\hun}(x_1, z) = \norm z -x_1\norm_{\un} = \inf\{t\geq 0 : z -x_1\leq t\un\}= \inf\{t\geq 0:  x_1\vee (x_2 + s\un) \leq x_1 + t\un\}$. If $t\geq 0$ and  $(x_2 + s\un) \leq x_1 + t\un$ then 
  $ x_1\vee (x_2 + s\un) \leq (x_1 + t\un)\vee x_1 =  (x_1 + t\un)$ which shows that 
  ${\sf D}_{\hun}(x_1, z) \leq \inf\{t\geq 0 : (x_2-x_1) \leqslant (t-s)\un\} = {\sf D}_{\hun}(x_1, x_2) + s = {\sf D}_{\hun}(x_1, x_2) - {\sf D}_{\hun}(z, x_2)$. $\Box$

\begin{prop}\label{geod1} For all $x_1, x_2\in E$ an for all $z\in\lb x_1, x_2\rb$
  \begin{equation*}
  {\sf D}_{\hun}(x_1, x_2) = {\sf D}_{\hun}(x_1, z) + {\sf D}_{\hun}(z, x_2)
  \end{equation*}
  \end{prop}
  
   \proof Let $z\in\lb x_1, x_2\rb$. By Lemma (\ref{eqsup}) either $w\in \lb x_1, x_1\vee x_2\rb$ or 
  $z\in \lb x_2, x_1\vee x_2\rb$. Without loss of generality assume that $z\in \lb x_1, x_1\vee x_2\rb$.\\ By Lemma \ref{splitbis} we have   
   ${\sf D}_{\hun}(x_1, x_1\vee x_2) = {\sf D}_{\hun}(x_1, z) + {\sf D}_{\hun}(z, x_1\vee x_2)$   
 and,    by Lemma \ref{split}, \\
   ${\sf D}_{\hun}(x_1, x_2) = {\sf D}_{\hun}(x_1, x_1\vee x_2) + {\sf D}_{\hun}(x_2, x_1\vee x_2)$. Therefore 
   \begin{equation*}
   {\sf D}_{\hun}(x_1, x_2) =  {\sf D}_{\hun}(x_1, z) + {\sf D}_{\hun}(z, x_1\vee x_2) + {\sf D}_{\hun}(x_2, x_1\vee x_2).
   \end{equation*} 
   From     Lemma \ref{eqsup} we  have ${\sf D}_{\hun}(z, x_1\vee x_2) = {\sf D}_{\hun}(z, z\vee x_2)$ and consequently  
   ${\sf D}_{\hun}(x_1, x_2) =  {\sf D}_{\hun}(x_1, z) + {\sf D}_{\hun}(z, z\vee x_2) + {\sf D}_{\hun}(x_2, z\vee x_2)$ which yields, by Lemma \ref{split},  ${\sf D}_{\hun}(x_1, x_2) =  {\sf D}_{\hun}(x_1, z) + {\sf D}_{\hun}(z, x_2)$. $\Box$

  \section{Geodesics}\label{geodesics}
  Given a metric space $(X, D)$, let us say that non empty subset $Z\subset X$ is  a geodesic, with respect to the metric $D$, if there exists an  onto map $\theta$ from a closed interval $[a, b]\subset\bbR$ to $Z$ such that, for all 
  $t_1, t_2\in [a, b]$, $D\big(\theta(t_1), \theta(t_2)\big) = \mid t_1 - t_2 \mid$ in which case we will say that 
  $\theta : [a, b]\to X$ is a parametrized geodesic from $\theta(a)$ to  $\theta(b)$.\\
    The metric space $(X, D)$ is a geodesic space if, for all pair $(x_1, x_2)\in X\times X$, there exists a parametrized geodesic $\theta : [a, b]\to X$ from $x_1$ to $x_2$.  
   A geodesic structure on  a geodesic metric space $(X, D)$ is a  family  $\Theta = \big(\theta_{(x_1, x_2)}, [a_{(x_1, x_2)}, b_{(x_1, x_2)}] \big)_{(x_1, x_2)\in X\times X}$ where 
  $\big(\theta_{(x_1, x_2)}, [a_{(x_1, x_2)}, b_{(x_1, x_2)}] \big)$ is a parametrized geodesic from $x_1$ to $x_2$  ; the geodesic $\theta_{(x_1, x_2)}\big([a_{(x_1, x_2)}, b_{(x_1, x_2)}\big)\subset X$ will be denoted by $\boldsymbol{\theta}(x_1, x_2)$.
  
   Given a geodesic structure  $\Theta$ on  a metric space $(X, D)$, we will say that a subset $C\subset X$ is geodesically convex (with respect to $\Theta$) if, for all $(x_1, x_2)\in X\times X$,  $\boldsymbol{\theta}(x_1, x_2)\subset C$; the empty set is convex by default. 
  
  If  $\theta : [a, b]\to X$ is a parametrized geodesic from $x_1$ to $x_2$ then, for all $s\in\bbR$, 
  $t\mapsto\theta(t-s)$ is a parametrized  geodesic $\tilde{\theta} $ on $[a+s, b + s]$ from $x_1$ to $x_2$ such that $\boldsymbol{\theta}(x_1, x_2) = \boldsymbol{\tilde{\theta}}(x_1, x_2)$ ; since  $D(x_1, x_2) =  b -  a$,  the map $\tilde{\theta} : [0, D(x_1, x_2)]\to X$ defined by 
  $\tilde{\theta}(t) = \theta(t+a)$ is a parametrized geodesic from from $x_1$ to $x_2$ with 
  $\boldsymbol{\theta}(x_1, x_2) = \boldsymbol{\tilde{\theta}}(x_1, x_2)$; let us call $\tilde{\theta} : [0, D(x_1, x_2)]\to X$ the standard parametrized geodesic  associated to the parametrized geodesic $\theta : [a, b]\to X$. 
  
  Given a  parametrized geodesic $\theta : [a, b]\to X$, from $x_1$ to $x_2$, the affinely parametrized geodesic \\$\hat{\theta} : [0, 1]\to X$ associated to  $\theta$  is   $\hat{\theta}(t) = \theta\big((1-t)a + tb\big) = 
  \theta\big(a + tD(x_1, x_2)\big)$.\\
   For all $t_1, t_2\in [0, 1]$, 
  $D\big(\hat{\theta}(t_1), \hat{\theta}(t_2) \big) = \mid t_1 - t_2\mid D(x_1, x_2)$.  
   
\bigskip  We show that $\big(E, {\sf D}_{\hun} \big)$ is a geodesic metric space, more precisely, there is a geodesic strucuture $\Gamma = \big(\gamma_{\un, (x_1, x_2)}, [\alpha_{\un, (x_1, x_2)}, \beta_{\un, (x_1, x_2)} ]\big)$ on $E$ for which the  geodesics are precisely the max-plus segments  and therefore, the max-plus convex sets are precisely the geodesicaly convex sets. To avoid cumbersome double subscripts we will write $\gamma_{\un}(x_1, x_2, t)$ for $\gamma_{\un (x_1, x_2)}(t)$ and similarly for $\alpha$ and $\beta$.  The explicit form of 
$\gamma_{\un}(x_1, x_2, t)$ is  given in Theorem \ref{geodesic}; for $E = \bbR^n$ and $\un = (1, \cdots, 1)$  this is  due to Stefan Gaubert, \cite{gauperso}.

 \begin{lem}\label{geod2}  Two arbitrary points $x_1, x_2\in E$ being given, any two other points of $\lb x_1, x_2\rb$ can be labelled    $z$ and  $z^\prime$ in such a way  one of the following two assertions holds:
 
 \begin{equation*}
 {\sf D}_{\hun}(z, z^\prime) =\begin{cases}
  {\sf D}_{\hun}(z, x_1\vee x_2) - {\sf D}_{\hun}(z^\prime, x_1\vee x_2) &\text{if } 
 x_i\leqslant z \leqslant z^\prime \leqslant x_1\vee x_2\\
  {\sf D}_{\hun}(z, x_1\vee x_2) + {\sf D}_{\hun}(z^\prime, x_1\vee x_2) &\text{if } x_1\leqslant z  \leqslant x_1\vee x_2 \text{ and }
    x_2\leqslant z^\prime  \leqslant x_1\vee x_2
  \end{cases}
  \end{equation*}

  \end{lem}
  
    \proof    From Lemma \ref{eqsup}, there are two cases to consider:  both $z$ and  $z^\prime$  are in $\lb x_i, x_1\vee x_2\rb$ for the same $i$,   one of $z$ and  $z^\prime$ is in $\lb x_1, x_1\vee x_2\rb$ and the other is in  $\lb x_2, x_1\vee x_2\rb$.\\
    Let us assume that $z$ and $z^\prime$ are both in  $\lb x_1, x_1\vee x_2\rb$ which, by     Lemma \ref{eqsup},  is a totally ordered set with respect the partial order of $E$.  If the two points  $z$ and $z^\prime$ are labeled such that   $x_1\leqslant z\leqslant z^\prime \leqslant x_1\vee x_2 $ then $z^\prime\in \lb z, x_1\vee x_2\rb$ and, from Proposition  \ref{geod1},   
    ${\sf D}_{\hun}(z, z^\prime) + {\sf D}_{\hun}(z^\prime, x_1\vee x_2) = {\sf D}_{\hun}(z, x_1\vee x_2)$

\medskip    For the second case, label the points such that  $z\in\lb x_1, x_1\vee x_2\rb$ and $z^\prime\in\lb x_2, x_1\vee x_2\rb$ 
   that is: \\ $x_1\leqslant z  \leqslant x_1\vee x_2$  and   
   $x_2\leqslant z^\prime  \leqslant x_1\vee x_2$ from which we have $ x_1\vee x_2 =  z\vee z^\prime$. Proposition  \ref{geod1} yields ${\sf D}_{\hun}(z, z^\prime) = {\sf D}_{\hun}(z, z\vee z^\prime) +  {\sf D}_{\hun}(z^\prime, z\vee z^\prime)$. $\Box$
   
   \begin{lem}\label{gama1to1} 
   If $x_1\leqslant x_2$ then the map $t\mapsto \gamma(t, x_1, x_2) =  x_1\vee (x_2 + t\un)$ is one to one and onto from the interval 
   $[\q_{\un}(x_1 - x_2), 0]$ to $\lb x_1, x_2\rb$. Furthermore, given $z\in\lb x_1, x_2\rb$,  the unique $s\in [\q_{\un}(x_1 - x_2), 0]$ such that $z = x_1\vee (x_2 + s\un)$ is $-{\sf D}_{\hun}(z, x_2)$. 
   \end{lem}
   
   \proof The points $x_1$ and $x_2$ being fixed, write $\gamma(t)$ for $\gamma(t, x_1, x_2)$.  If $z\in \lb x_1, x_2\rb$ then either $z = (x_1 + t\un)\vee x_2$ or $z = x_1\vee (x_2 + t\un)$ with $t\leq 0$ ; in the first case $z = x_2 = x_1\vee x_2 = x_1\vee (x_2 + 0\un)$. This shows that $\gamma$ is onto. 
   
   \medskip To see that $\gamma$ is one to one on $[\q_{\un}(x_1 - x_2), 0]$ let $z = x_1\vee (x_2 + s\un)$, 
   $z^\prime = x_1\vee (x_2 + t\un)$ with $\q_{\un}(x_1 - x_2)\leq s \leq t\leq 0$ and assume that $z = z^\prime$. We show that $s = t$. 
   
   \medskip First, $x_2 + t\un \leqslant  x_1\vee (x_2 + t\un) = z^\prime = z$ from which we have 
   $t\un \leqslant z - x_2$ and therefore,   from the definition of $\q_{\un}$, 
   \begin{equation}\label{tqox2}
   t\leq \q_{\un}(z - x_2)
   \end{equation} 
   
   Now, $z - x_2  = \big[x_1\vee (x_2 + s\un)\big] - x_2 =  (x_1-x_2)\vee \big((x_2 + s\un) - x_2 \big) = (x_1-x_2)\vee (s\un)$. 
   
   \medskip To complete the proof, we identify $E$ with a Riesz subspace of the space $\mathcal{C}(\Omega)$ containing the  constant function $\boldsymbol{1}$,  where $\Omega$ is a compact Hausdorff topological space, and we set 
   $\un = \boldsymbol{1}$. 
   
    \begin{equation}
  \displaystyle{ \q_{\un}(x_1 - x_2) = \inf_{\omega\in\Omega}\big(x_1(\omega) - x_2(\omega)\big)\, \leq\,  s}
   \end{equation}
   
   and therefore, 
   \begin{equation}\label{omegx22}
   \q_{\un}\big( (x_1-x_2)\vee (s\boldsymbol{1})\big) = \inf_{\omega\in\Omega}\max\left\{x_1(\omega) - x_2(\omega), s\right\} \leq s
   \end{equation}

   We have shown that $\q_{\un}(z - x_2) \leq s$ which, with (\ref{tqox2}), gives $s =t$. 
   
   \bigskip The last part follows from Lemma \ref{wxypu}.$\Box$

 \bigskip For $x_1, x_2\in E$ let $\alpha_{\un}(x_1, x_2) = \min\{0, \q_{\un}(x_1 - x_2)\}$ and $\beta_{\un}(x_1, x_2) = \max\{0, \p_{\un}(x_1-x_2)\}$ Then $ \beta_{\un}(x_1, x_2)= -\alpha_{\un}(x_2, x_1)$,  
$\alpha_{\un}(x_1, x_2)\leq\beta_{\un}(x_1, x_2)$ and ${\sf D}_{\hun}(x_1, x_2) = \beta_{\un}(x_1, x_2) - \alpha_{\un}(x_1, x_2)$.

 \begin{lem}\label{mainprop}  If $x_1$ and $x_2$ are  two elements  of $E$ and $s, t$ are real numbers such that  $\alpha_{\un}(x_1, x_2)\leq t \leq 0 \leq s \leq \beta_{\un}(x_1, x_2)$ then,   for  
    $z = x_1\vee (x_2 + t\un)\in\lb x_1, x_1\vee x_2\rb$ and $z^\prime = (x_1 -  s\un) \vee x_2\in\lb x_2, x_1\vee x_2\rb$ one has 
  \begin{equation} \begin{cases}
 (1)\quad  {\sf D}_{\hun}(z, x_1\vee x_2) = -t&\\ 
  &\text{and}\\
 (2)\quad  {\sf D}_{\hun}(z^\prime, x_1\vee x_2) = s&
  \end{cases}
  \end{equation}
  \end{lem}
 
 \proof 
 We have $z =  x_1\vee (x_2 + t\un) \leqslant  x_1\vee\big( (x_1\vee x_2) + t\un\big) = x_1\vee [(x_1 + t\un)\vee (x_2+t\un)] = 
 [x_1\vee (x_1 + t\un)]\vee (x_2+t\un) = x_1\vee (x_2+t\un)$.

\medskip\noindent In conclusion,  $z =  x_1\vee\big( (x_1\vee x_2) + t\un\big)$.

\medskip Let us see that 
 $\min\{0, \q_{\un}(x_1 - x_2)\} =  \q_{\un}(x_1 - x_1\vee x_2)$. 
 
 \medskip\noindent
 By definition $\q_{\un}(x_1 - x_2) = \sup\{t : t\un \leqslant x_1 - x_2\}$. If $t\leq 0$ and $t\un \leqslant x_1 - x_2$ then 
 $x_2\leqslant x_1 - t\un$ and $x_1\vee x_2 \leqslant(x_1 -t\un)\vee x_1 = x_1- t\un $ which shows that 
 $\min\{0, \q_{\un}(x_1 - x_2)\}\leq \q_{\un}(x_1 - x_1\vee x_2)$.
 
 \medskip\noindent
 Reciprocally, if $t\leq  \q_{\un}(x_1 - x_1\vee x_2)$ then $t\leq 0$ and $t\un \leqslant (x_1 - x_1\vee x_2) \leqslant x_1 - x_2$ and therefore, \\$t\leq \min\{0, \q_{\un}(x_1 - x_2)\}$.
 
  \medskip The conclusion follows from Lemma \ref{gama1to1}  and $ \beta_{\un}(x_1, x_2)= -\alpha_{\un}(x_2, x_1)$. $\Box$
 
\begin{lem}\label{lemgamma}
 The map $\eta_{\un}$ defined on $E\times E\times\bbR$ by  
 \begin{equation}\label{defgamma}\eta_{\un}(x_1, x_2, t) = 
 \begin{cases}
 x_1\vee (x_2 + t\un)&\text{if }   t\leq 0\\
 (x_1 -  t\un) \vee x_2&\text{if } 0\leq t
 \end{cases}
 \end{equation}
is uniformly continuous with respect to the product topology induced by the metric $ {\sf D}_{\hun}$ on $E$. Furthermore, 
\begin{equation}\label{propgamma}
\begin{cases}
\eta_{\un}(x_1, x_2, t) = x_1 \text{ if } t\leq \alpha_{\un}(x_1, x_2) \text{ and } 
\eta_{\un}(x_1, x_2, t) = x_2 \text{ if } \beta_{\un}(x_1, x_2) \leq t \\
\eta_{\un}\left(x_1, x_2, \bbR\right) = \lb x_1, x_2\rb = \eta_{\un}\left(x_1, x_2,  [\alpha_{\un}(x_1, x_2), \beta_{\un}(x_1, x_2)]\right)
\end{cases}
\end{equation} 
 \end{lem}
 
 \proof Since $\norm \cdot \norm_{\un}$ is a Riesz norm, the lattice operations are uniformly continuous with respect to  $\norm \cdot \norm_{\un}$ from which it follows that $\eta_{\un}$ is uniformly continuous with respect to $\norm \cdot \norm_{\un}$ and therefore with respect $ {\sf D}_{\hun}$, since the norms $\norm \cdot \norm_{\un}$ and $\norm \cdot \norm_{\hun}$ are equivalent. 
 The verification of (\ref{propgamma}) is straightforward. $\Box$

 \begin{thm}\label{geodesic}
  For all $x_1, x_2\in E$ and $t\in\bbR$ let $\gamma_{\un}$ be the restriction of $\eta_{\un}$ to $ [\alpha_{\un}(x_1, x_2),  \beta_{\un}(x_1, x_2)]$ that is, 
 \begin{equation}\label{defgamma}\gamma_{\un}(x_1, x_2, t) = 
 \begin{cases}
 x_1\vee (x_2 + t\un)&\text{if }  \alpha_{\un}(x_1, x_2)\leq t\leq 0\\
 (x_1 -  t\un) \vee x_2&\text{if } 0\leq t\leq  \beta_{\un}(x_1, x_2)
 \end{cases}
 \end{equation}
Then,  for all $x_1, x_2\in E$ and  for all  $t_1, t_2\in [\alpha_{\un}(x_1, x_2),  \beta_{\un}(x_1, x_2)]$
 \begin{equation} {\sf D}_{\hun}\big(\gamma_{\un}(x_1, x_2, t_1), \gamma_{\un}(x_1, x_2, t_2)\big) = \vert t_1 - t_2\vert \end{equation}
 In other words, $\Gamma_{\un} = \big(\gamma_{\un}(x_1, x_2, -), [\alpha(x_1, x_2), \beta(x_1, x_2)] \big)_{(x_1, x_2)\in X\times X}$ is  a  geodesic structure on $\big(E, {\sf D}_{\un}\big)$ for which the geodesics are precisely the max-plus segments.
  \end{thm}
  
 \proof  Let $z = \gamma_{\un}(x_1, x_2, t)$ and $z^\prime = \gamma_{\un}(x_1, x_2, s)$. 
   From Lemma \ref{eqsup}, there are three  cases to consider:
   
   $$\begin{cases}
   (1)\, \,  z, z^\prime\in \lb x_1, x_1\vee x_2\rb\\
   (2)\, \,  z, z^\prime\in \lb x_2, x_1\vee x_2\rb\\
   (3)\, \,  z\in \lb x_1, x_1\vee x_2\rb \hbox{ and }  z^\prime\in \lb x_2, x_1\vee x_2\rb.
   \end{cases}$$
   
   In all cases, the conclusion follows from Lemma  
   \ref{mainprop} and
Lemma \ref{geod2}. \hfill$\Box$

\begin{coro}\label{m+conv=geoconv}
 A subset $C$ of $E$ is max-plus convex if and only if it is a geodesically convex set with respect to the geodesic structure $\Gamma_{\un}$. 
 \end{coro}

 \bigskip Let $\hat{\gamma}_{\un} : E\times E\times [0, 1]\to E$ be the  affinely parametrized geodesic associated to $\gamma_{\un}$ that is,  
$$\hat{\gamma}_{\un}(x_1, x_2, t) = \gamma_{\un}\big(x_1, x_2, (1-t)\alpha_{\un}(x_1, x_2) + t\beta_{\un}(x_1, x_2)\big) = 
 \gamma_{\un}\big(x_1, x_2, \alpha_{\un}(x_1, x_2) + t{\sf D}_{\un}(x_1, x_2)\big)$$
then 
\begin{equation}\label{normgeo}
{\sf D}_{\hun}\big(\hat{\gamma}_{\un}(x_1, x_2, s), \hat{\gamma}_{\un}(x_1, x_2, t)\big) = \mid s-t \mid\, {\sf D}_{\hun}(x_1, x_2)
\end{equation}

\bigskip If one defines the midpoint map $\mu_{\un} : E\times E\to E$ by $\mu_{\un}(x_1, x_2) = \hat{\gamma}_{\un}(x_1, x_2, 1/2)$ then: 
\begin{equation}\label{propmidptmp}
\forall (x_1, x_2)\in E\times E\quad
\begin{cases}
(1)& \mu_{\un}(x_1, x_2) = \mu_{\un}(x_2, x_1)\\
&\\
(2)&\displaystyle{{\sf D}\left(x_1, \mu_{\un}(x_1, x_2)\right) = {\sf D}\left(x_2, \mu_{\un}(x_1, x_2)\right) = \frac{1}{2}{\sf D}\left(x_1, \mu_{\un}(x_1, x_2)\right)}
\end{cases}
\end{equation}

A closed subset $C$ on a topological vector space is convex (in the usual affine structure) if and only if, for all $x_1, x_2\in C$, $\displaystyle{\frac{1}{2}x_1 + \frac{1}{2}x_2}\in C$. 

\medskip Proposition \ref{carmx+mid} below characterizes closed max-plus convex 
sets as  semilattices containing the (max-plus) midpoints of any pair of their points. 

\medskip First, notice that, if $x_1\leqslant x_2$ then $-\alpha_{\un}(x_1, x_2) = {\sf D}_{\un}(x_1, x_2) = \norm x_2 - x_1\norm_{\un} = \beta_{\un}(x_2, x_1)$ and\\
 $x_1\leqslant\mu_{\un}(x_1, x_2) \leqslant x_2 $ since $\mu_{\un}(x_1, x_2) = x_1\vee \big(x_2 + (\alpha_{\un}(x_1, x_2)/2)\un\big)$.
 
 \begin{lem}\label{unifcontgammahat} With respect to the metric topology associated to ${\sf D}_{\un}$ on $E$, the affinely parametrized geodesic $\hat{\gamma}_{\un} : E\times E\times [0, 1]\to E$ is continuous on $ E\times E\times [0, 1]$. 
\end{lem}

\proof  The restriction of $\eta_{\un}\big(x_1, x_2, \alpha_{\un}(x_1, x_2) + t{\sf D}_{\un}(x_1, x_2)\big)$  to $E\times E\times [0, 1]$ is $\hat{\gamma}_{\un}(x_1, x_2, t)$ and, by Lemma \ref{lemgamma},  $\eta_{\un} : E\times E\times\bbR\to E$ is  continuous.  Showing that $(x_1, x_2)\mapsto \alpha_{\un}(x_1, x_2)$ is  continuous on $E\times E$ will complete the proof. We show that $(x_1, x_2)\mapsto \beta_{\un}(x_1, x_2)$ is uniformly continuous on $E\times E$ with respect to the topology of the norm $\norm \cdot \norm_{\un}$, which is a Riesz norm on $E$, which implies that $x\mapsto x^+$ is  uniformly continuous on $E$, and therefore $x\mapsto \norm x^+ \norm_{\un}$ is  uniformly continuous. The conclusion follows from 
$\norm x^+ \norm_{\un} = \p_{\un}(x^+) = \max\{0, \p_{\un}(x)\} = \beta_{\un}(x, 0)$ and 
$\beta_{\un}(x_1, x_2) = \beta_{\un}(x_1-x_2, 0)$. $\Box$

\begin{coro}\label{contmid}
The midpoint map $\mu_{\un} : E\times E\to E$ is continuous.
\end{coro}

\begin{prop}\label{carmx+mid}For all non empty closed subset $C$ of $E$  the first and the last of the three assertions  below  are equivalent; if $C$ is complete then the three assertions are equivalent.  

\medskip\noindent{\bf(A)} $C$ is max-plus convex.  

\medskip\noindent{\bf(B)} $
\forall x_0, x_1\in C\quad \mu_{\un}(x_0, x_1)\in C.$

\medskip\noindent{\bf(C)}
$\forall x_0, x_1\in C\quad 
\begin{cases} 
(1)& x_0 \vee x_1\in C \text{ and }\\
(2)& \text{ if } x_0\leqslant x_1 \text{ then }\mu_{\un}(x_0, x_1)\in C.
\end{cases}$
\end{prop}

\proof Clearly, $(A)$ implies $(B)$ and $(C)$. We show that $(C)$ implies $(A)$.

We have to see that, for all $x_0, x_1\in C$, $\lb x_0, x_1\rb\subset C$. From Lemma \ref{mainconvlem2} and hypothesis $(1)$, it is sufficient to show that $\lb x_0, x_1\rb\subset C$ whenever $x_0, x_1\in C$ and $x_0\leqslant x_1$.\\
Let $D_n = \{(k/2^n) : 0\leq k \leq n\}$. If $x_0, x_1\in C$ and $x_0\leqslant x_1$ then, by the second hypothesis, $x_{1/2} = \mu_{\un}(x_0, x_1)\in C$ and $x_0\leqslant x_{1/2}\leqslant x_1$.\\
Let $n\geq 1$ and assume that, for each $m\leq n$ we have a sequence of points $S_m = \{x_{m, t} : t\in D_m\}\subset \lb x_0, x_1\rb$ such that: \\
$(a)$ If $t\in D_{m-1}$ then $x_{m-1, t} = x_{m, t}$\\
$(b)$ If $t, t^\prime\in D_{m}$ and $t\leq t^\prime$ then $x_{m, t}\leqslant x_{m, t^\prime}$\\
$(c)$ If ${\sf D}_{\hun}(x_{n, k/2^n}, x_{n, (k+1)/2^n}) = 2^{-n}{\sf D}_{\hun}(x_0, x_1)$.\\
To construct $S_{n+1}$ such that $(a), (b)$ and $(c)$ hold consider $t = (k/2^{n+1})$ and $t\not\in D_n$; there  are then  $t^{\prime}, t^{\prime\prime}\in D_n$ such that $t = (t^{\prime} + t^{\prime\prime})/2$; let $x_{n+1, t}= \mu_{\un}(x_{n, t^{\prime}}, 
x_{n, t^{\prime\prime}})$. \\
Now, let $S = \{x_t : t\in D\} = \cup_{n}S_n$ where $D\subset [0, 1]$ is the set of dyadic numbers. For all given $n$, $S_n$ is a linearly ordered subset of 
$\lb x_0, x_1\rb$. Take $x\in \lb x_0, x_1\rb$; if $x = x_0$ then ${\sf D}_{\hun}(x, x_{n, 1/2^n}) = (1/2^n){\sf D}_{\hun}(x_0, x_1)$; if $x = x_1$ then ${\sf D}_{\hun}(x, x_{n, 2^{n-1}/2^n}) = (1/2^n){\sf D}_{\hun}(x_0, x_1)$. If $x$ is neither $x_0$ nor $x_1$, nor a point of $S_n$,   there are  two points $x_{n, k/2^n}$ and  $x_{n, (k+1)/2^n}$ such that $x_{n, k/2^n} \leqslant x \leqslant x_{n, (k+1)/2^n}$, since, by Lemma \ref{eqsup}, $\lb x_0, x_1\rb$ is linearly ordered by the restriction of the partial of $E$.  We have 
$\lb x_0, x_1\rb = \lb x_0, x_{n, k/2^n}\rb \cup \lb x_{n, k/2^n},  x_{n, (k+1)/2^n}\rb \cup \lb x_{n, (k+1)/2^n}, x_1\rb$ and since $x\not\in S_n$, $x$ belongs to only one these max-plus segments, which is $ \lb x_{n, k/2^n},  x_{n, (k+1)/2^n}\rb$. From $(1/2^n) = {\sf D}_{\hun}( x_{n, k/2^n},  x_{n, (k+1)/2^n}) =  {\sf D}_{\hun}( x_{n, k/2^n},  x) + {\sf D}_{\hun}( x,  x_{n, (k+1)/2^n})$ we have ${\sf D}_{\hun}( x_{n, k/2^n},  x) \leq (1/2^n)$. \\
We have shown that $S$ is a dense subset of $\lb x_0, x_1\rb$ and that $S\subset C$. Since $C$ is closed we have 
$\lb x_0, x_1\rb\subset C$. 

\medskip Assuming that $C$ is complete - with respect to ${\sf D}_{\hun}$ - we show that $(B)$ implies $(A)$.\\
Since $(C, {\sf D}_{\hun}, \mu_{\un})$ is a complete midpoint space in the sense of \cite{midpt} there exists an affinely parametrized geodesic $\varphi : C\times C \times [0, 1]\to C$ such, for all $(x_0, x_1)\in C\times C$,  $\mu_{\un}(x_0, x_1) = \varphi(x_0, x_1, 1/2)$ which is obtained by dyadic approximation starting from $\mu_{\un}(x_1, x_2)$, Lemma 3.0.1 and its proof in \cite{midpt}; the restriction of $\hat{\gamma}_{\un}$ to $C\times C \times [0, 1]$ is also an affinely parametrized geodesic. Furthermore, for all  $(x_0, x_1)\in C\times C$, $ \hat{\gamma}_{\un}(x_0, x_1, 1/2)=\mu_{\un}(x_0, x_1) = \varphi(x_0, x_1, 1/2)$. By dyadic approximation we have, 
$(x_0, x_1)\in C\times C\times [0, 1]$, $\hat{\gamma}_{\un}(x_0, x_1, t)= \varphi(x_0, x_1, t)\in C $ that is 
$\lb x_0, x_1\rb\subset C$.
$\Box$

\section{On the topology of max-plus convex sets}\label{topmax}

In this section, the topology on a given  Riesz space $E$ is the metric topology ${\sf D}_{\hun}$ associated to a given unit $\un$ of $E$. Either of the first proposition or the first lemma of this section shows that in  the metric space  $\big(E,  {\sf D}_{\hun}\big)$ open and closed balls are max-plus convex, and therefore geodesically convex ; consequently,  the topology is locally 
max-plus convex.  A max-plus convex set is an absolute retract from which we have the max-plus version of the classical Kakutani Fixed Point Theorem. There are also max-plus versions of some classical continuous (approximate) selection theorems for upper semicontinuous multivalued maps (here simply called ``maps'').

\begin{lem}\label{hulocmxconv}
Balls, open or closed  with respect to $\norm \cdot \norm_{\hun}$,  are max-plus convex. 
\end{lem}

\proof {\bf (A)} 
Let $B = B_{\hun}(0, \delta) = \big\{x\in E :  \norm x \norm_{\hun} < \delta \big\}$; from  Theorem \ref{hunrienrm} we have 
$x_1\vee x_2\in B$ if $\{x_1, x_2\}\subset B$.

\medskip Since $\lb x_1, x_2\rb = \lb x_1, x_1\vee x_2\rb\, \cup\, \lb x_1, x_1\vee x_2\rb$ it is sufficient to prove that $\lb x_1, x_2\rb\subset B$ whenever $x_1, x_2\in B$ and $x_1$ and $x_2$ are comparable.

\medskip{\bf(1)} Assume that $x_1\leqslant x_2$. 

\medskip If $t\leq 0$ then $x_1 + t_1\un \leq x_2$  from which we have 
$( x_1 + t_1\un)\vee x_2 = x_2\in B$. 

\medskip{\bf(2)} 
Assume that $x_2\leqslant x_1$. 

\medskip We have to see that, if $t \leq 0$ then  $(x_1 + t\un)\vee x_2\in B$  or, equivalently, that 
$(x_1 - s\un)\vee x_2\in B$ if $s\geq 0$.

From $s\geq 0$ and  $x_2\leqslant x_1$ we have  $(x_1 - s\un)\vee x_2\leq x_1\vee x_2 = x_1$ and therefore 
$\p_{\un}\big( (x_1 - s\un)\vee x_2\big)\leq \p_{\un}(x_1)$ and 
\begin{equation}\label{s121}
\max\{0,\p_{\un}\big( (x_1 - s\un)\vee x_2\big)\}\leq \max\{0, \p_{\un}(x_1)\}
\end{equation}

\medskip From $(s\un -x_1)\wedge (-x_2)\wedge (-x_1) = (s\un -x_1)\wedge (-x_1)\leq (-x_1)$ and $-[(x_1 - s\un)\vee x_2] =(s\un -x_1)\wedge (-x_2)$ we have $-[(x_1 - s\un)\vee x_2] \leqslant (-x_1)$ 
 from which, $\p_{\un}\big(-[(x_1 - s\un)\vee x_2]\big)\leq \p_{\un}(-x_1)$ and 
\begin{equation}\label{s122}\max\{0, \p_{\un}\big(-[(x_1 - s\un)\vee x_2]\big)\}\leq \max\{0, \p_{\un}(-x_1)\}
\end{equation} 

\medskip Adding the inequalities from (\ref{s121}) and (\ref{s122})  gives $\norm (x_1 - s\un)\vee x_2] \norm_{\hun} \leq \norm x_1\norm_{\hun}$.  

\medskip We have shown that $B$ is max-plus convex; the same procedure shows that closed balls centered at $0$ are max-plus convex;

 \bigskip {\bf(B)} We show that arbitrary balls are max-plus convex. 
 Let $B = \{x\in E : {\sf D}_{\hun}(x_0, x) \leq\delta\}$ be a ball centered at $x_0$ and let $x_1$ and $x_2$ be two points of $B$; from $x_i -x_0\in (B - x_0)$ and from the first part of the proof, we have, for all $t\leq 0$, $\big[(x_1 - x_0) + t\un\big]\vee( x_2 -x_0) \in (B - x_0)$ and therefore, $(x_1  + t\un )\vee x_2  \in B $.
$\Box$

\begin{prop}\label{quasiconmet} 
For all $(y, x_1, x_2)\in E\times E$ and for all $x\in\lb x_1, x_2\rb$
\begin{equation}\label{quasimid}
{\sf D}_{\hun}\big(y, x\big) \leq \max\left\{{\sf D}_{\hun}(y, x_1), {\sf D}_{\hun}(y, x_2) \right\}
\end{equation}
and, more generally, for all non empty  set $S\subset E$ and for all $y\in E$
\begin{equation}\label{quasimidfin}
\forall x\in \lb S\rb\quad {\sf D}_{\hun}\big(y, x\big) \leq \max_{z\in S} {\sf D}_{\hun}\big(y, z\big) 
\end{equation}
\end{prop}

\proof Let $r =  \max\left\{{\sf D}_{\hun}(y, x_1), {\sf D}_{\hun}(y, x_2) \right\}$ and let $B$ be the closed ball - with respect to 
${\sf D}_{\hun}$ - of radius $r$ centered at $y$; from $x_i\in B$ and Lemma \ref{hulocmxconv} we have 
$\lb x_1, x_2\rb\subset B$ which establishes (\ref{quasimid}). 

\medskip\noindent Assume $S = \{ x_1, \cdots, x_{m}\} $ with $m > 2$ and let $S_1 = \{ x_1, \cdots, x_{m-1}\}$ and proceed by induction: if 
$x\in \lb S_1\rb$ then ${\sf D}_{\hun}\big(y, x\big) \leq \max_{1\leq i\leq m-1} {\sf D}_{\hun}\big(y, x_i\big)$; if $x\in\lb S\rb$ then, by Lemma \ref{recur}, $x\in\lb x^\prime, x_m\rb$ with $x^\prime\in\lb S_1\rb$ and, by (\ref{quasimid}), 
 ${\sf D}_{\hun}\big(y, x\big) \leq \max\{{\sf D}_{\hun}\big(y, x^\prime\big), {\sf D}_{\hun}\big(y, x_m\big)\}$. We have shown that (\ref{quasimidfin}) holds for finite sets.
 
 \medskip\noindent For the general case, from 
  $\displaystyle{\lb S \rb = \bigcup_{A\in\langle S \rangle}\lb A \rb}$, $(3)$ of Lemma \ref{mainconvlem1},  we have $x\in\lb A\rb$ for some finite set 
  $A\subset S$ from which the conclusion follows. 
 $\Box$
 
 \begin{coro}\label{diam}
 The diameter with respect to the metric ${\sf D}_{\hun}$ of $\lb S\rb$ is $\max\left\{{\sf D}_{\hun}(x, y) : (x, y)\in S\times S\right\}$ which is the diameter of $S$. 
 \end{coro}
 
 \bigskip Lemma \ref{hulocmxconv} and Proposition \ref{quasiconmet} are equivalent, they say that the metric  
 ${\sf D}_{\hun}$ is (max-plus) quasiconvex. One can obviously define a max-plus quasiconvex function $f : C\to\bbR$, where $C$ is a max-plus convex subset of $E$, by the property $\max_{x\in\lb x_1, x_2\rb}f(x) = \max\{f(x_1), f(x_2)\}$ and show that, for such functions, $\max_{x\in\lb S\rb}f(x) = \max_{x\in S}f(x)$.

\bigskip Given two subsets $S_1, S_2$ of $E$ let $S_1\vee S_2 = \{v_1\vee v_2 : (v_1, v_2)\in S_1\times S_2\}$ and let 
$S_i + t\un = \{v_i +  t\un : v_i\in S_i\}$.

\begin{lem}\label{bVb}
If $(E, \norm \cdot \norm_{\un})$ is  complete  then, 
\begin{equation*}
\forall\,  x_1, x_2\in E\quad B_{\un}(x_1, \delta) \vee B_{\un}(x_2 , \delta) = B_{\un}(x_1 \vee x_2, \delta)
\end{equation*}
where $B_{\un}(x, \delta)$ is the open ball with respect to $ \norm \cdot \norm_{\un}$ of radius $\delta$ centered at 
$x$.
\end{lem}

\proof We can assume that $(E, \norm \cdot \norm_{\un})$ is $\mathcal{C}(\Omega)$, for some compact topological space $\Omega$, and that $\un : \Omega\to\bbR$ is the constant function $1$; $\norm \cdot \norm_{\un}$  is then  the  sup-norm $\norm x \norm_{\infty} = \max_{\omega\in\Omega}\vert x(\omega)\vert$; we simply write $B(x, \delta)$ for the open ball with respect to the sup-norm.  If $\norm x_i - y_i\norm_{\infty} < \delta$ then, for all $\omega\in\Omega$, $x_i(\omega) - \delta < y_i(\omega) < x_i(\omega) + \delta$  and therefore 
$\max\{x_1(\omega), x_2(\omega)\} - \delta < \max\{y_1(\omega), y_2(\omega)\} < \max\{x_1(\omega), x_2(\omega)\} + \delta$ which shows that  $B(x_1, \delta) \vee B(x_2 , \delta) \subseteq B(x_1 \vee x_2, \delta)$.\\
To prove the other inclusion let $\Omega_1 = \{\omega\in\Omega : x_2(\omega)\leq x_1(\omega)\}$ and similarly for $\Omega_2$.\\
If $y\in B_{\un}(x_1 \vee x_2, \delta)$ then 
\begin{equation*}
\forall \omega\in\Omega_i\quad  - \delta <  y(\omega) - x_i(\omega) < \delta.
\end{equation*}
By the Tietze-Urysohn's Theorem, there exists a continuous function $v_i : \Omega\to ]-\delta, \delta[$ such that, 
\begin{equation*}
\forall \omega\in\Omega_i\quad  v_i(\omega) =  y(\omega) - x_i(\omega).
\end{equation*}
Let $z_i = v_i + x_i$ and notice that $\begin{cases}
\forall \omega\in\Omega_i\quad z_i(\omega) = y(\omega)\\
\forall \omega\in\Omega\quad x_i(\omega)-\delta < z_i(\omega) < x_i(\omega)+\delta
\end{cases}$\\
from which we have $z_i\in B(x_i, \delta)$ and $y\leq \max\{z_1, z_2\}$. Let $y_i = \min\{y, z_i\}$ then 
$y_i\in B(x_i, \delta)$ and $y = \max\{y_1, y_2\}$. $\Box$

\begin{prop}\label{intmx+}
If $\norm \cdot \norm_{\un}$ is a complete norm on $E$ then the interior of a max-plus convex subset of $E$ is max-plus convex.
\end{prop}

\proof Let $C\subseteq E$ be a max-plus convex set whose interior $\stackrel{\circ}{C}$ is not empty; take $x_1, x_2\in \stackrel{\circ}{C}$ and $\delta > 0$ such that $B_{\un}(x_i, \delta)\subset C$. For all $t\leq 0$, 
$B_{\un}(x_1, \delta)\vee\big(B_{\un}(x_2, \delta) + t\un\big)\subset C$ and $B_{\un}(x_2, \delta) + t\un = B_{\un}(x_2 + t\un, \delta)$ therefore, $B_{\un}(x_1, \delta)\vee  B_{\un}(x_2 + t\un, \delta)\subset C$ which, from Lemma \ref{bVb}, gives 
$B_{\un}\big(x_1\vee(x_2 + t\un), \delta)\subseteq C$. $\Box$

\bigskip Given $\delta > 0$ and a non empty set $S\subseteq E$ let 
$U_{\delta}(S) = \{x\in E : \inf_{y\in S}\norm x-y\norm_{\un} < \delta\}$ and \\$V_{\delta}(S) = \{x\in E : \inf_{y\in S}\norm x-y\norm_{\hun} < \delta\}$

\begin{lem}\label{unifconv} For all non empty subset $E$ of $E$ one has $V_{\delta}(S)\subset U_{\delta}(S) \subset V_{2\delta}(S)$ and, for all  max-plus convex set $C\subset E$ and for all $\delta > 0$, $U_{\delta}(C)$ is max-plus convex.  As a consequence, the closure of a max-plus convex set is max-plus convex. 
\end{lem}

\proof The first  part follows from $\norm \cdot \norm_{\un} \leq \norm \cdot \norm_{\hun} \leq 2 \norm \cdot \norm_{\un}$.\\
Given a non empty max-plus set $C\subset  E$, take $y_i\in  U_{\delta}(C)$ and $x_i\in C$, $i=1, 2$,  and $\eta > 0$ such that  $\norm x_i-y_i\norm_{\un} \leq \eta < \delta$.  From   
$\bvert  x_i-y_i \bvert \leqslant \eta \un$ we have $x_i - \eta\un \leqslant y_i \leqslant x_i + \eta\un$. For all $t\in\bbR$, 
$x_2 + t\un - \eta\un \leqslant y_2 + t\un \leqslant x_2 + t\un + \eta\un$.  
From $x_2 + t\un - \eta\un \leqslant y_2 + t\un$ and 
$x_1-\eta\un \leqslant y_1 $ we have 
$(x_1 -\eta\un)\vee (x_2 + t\un - \eta\un)\leqslant y_1\vee (y_2 + t\un)$ that is 
$x_1\vee(x_2 + t\un) - \eta\un \leqslant y_1\vee (y_2 + t\un)$. Similarly, from $(y_2 + t\un) \leqslant (x_2 + t\un + \eta\un)$ and $y_1\leqslant x_1 + \eta\un$ we have 
$y_1\vee (y_2 + t\un) \leqslant x_1\vee (x_2 + t\un) + \eta\un$. We have shown that 
$\bvert y_1\vee (y_2 + t\un) -  x_1\vee (x_2 + t\un)  \bvert \leqslant \eta\un$. If $t\leq 0$  then 
 $x_1\vee (x_2 + t\un)\in C$ and $\norm y_1\vee (y_2 + t\un) -  x_1\vee (x_2 + t\un)  \norm_{\un} \leq \eta$ which shows that $y_1\vee (y_2 + t\un)\in   U_{\delta}(C)$.\\
 The last part follows from $\overline{C} = \cap_{\delta > 0}U_{\delta}(C)$. $\Box$
 
 \begin{prop}\label{toppolymax}$\phantom{i}$\\ 
$(1)$ If $x_1\neq x_2$ then $\lb x_1, x_2\rb$ is a topological arc in the metric space 
$\big(E, {\sf D}_{\un}\big)$.\\
$(2)$ For all finite sets $\{x_1, \cdots, x_m\}\subset E$, $\lb x_1, \cdots, x_m\rb$ is a compact subset of the metric space 
$\big(E, {\sf D}_{\un}\big)$. 
 \end{prop}
 
 \proof $(1)$ If $x_1\neq x_2$ then $\alpha_{\un} < \beta_{\un}$ and $\gamma_{\un}: [\alpha_{\un}(x_1, x_2), \beta_{\un}(x_1, x_2)]\to \lb x_1, x_2\rb$ is a homeomorphism.
 
 \bigskip\noindent $(2)$ By induction on $m$. We can assume that $m\geq 3$. By the induction hypothesis 
 $K = \lb x_1, \cdots, x_{m-1}\rb$ is compact. Take $t_1, t_2\in\bbR$ such that, for all $x\in K$, 
 $t_1\un \leqslant x - x_m \leqslant t_2\un$, Lemma \ref{contlatop}. We have, for all $x\in K$, 
 $t_1\leq \q_{\un}(x-x_m)$ and $\p_{\un}(x-x_m)\leq t_2$; without loss of generality, we can assume that $t_1 < 0 < t_2$ which yields, for all $x\in K$, $t_1\leq\alpha(x, x_m) \leq \beta(x, x_m) \leq t_2 $.\\ 
 On $K\times [t_1, t_2]$ define $\gamma_\star$ by $\gamma_\star(x, t) = \gamma_{\un}(x, x_m, t)$; from 
 $[\alpha(x, x_m), \beta(x, x_m)] \subseteq [t_1, t_2]$ we have, for all $x\in K$, $\gamma_\star\big(x, [t_1, t_2]\big) = \lb x, x_m\rb$ and, from Lemma \ref{recur}, we have $\gamma_\star\Big(K\times [t_1, t_2]\Big) = \lb x_1, \cdots, x_m\rb$. \\
 To complete the proof, notice that by the definition of $\gamma_{\un}$, (\ref{defgamma}), and Lemma \ref{contlatop}, 
 $\gamma_\star$ is continuous. $\Box$
 
 \begin{lem}\label{homotriv} Non empty  max-plus convex susbets of  $E$  are  contractible. 
\end{lem}

\proof By Lemma \ref{unifcontgammahat}, $\hat{\gamma} : E\times E\times [0, 1]\to E$ is continuous ; if  $C$ be a max-plus subset of $E$ then $\hat{\gamma}\big(C\times C\times [0, 1]\big) = C$ futhermore, for all $(x, y)\in C\times C$, $\hat{\gamma}(x, y, 0) = x$ and  $\hat{\gamma}(x, y, 1) = y$. $\Box$

\begin{prop}\label{hullcomp}
If $C\subset E$ is a max-plus convex subset of $E$ which is complete with respect to ${\sf D}_{\hun}$ then, for all compact subsets $K\subset C$, $\lb K \rb$ is compact. 
\end{prop}

\proof From Proposition \ref{toppolymax} and Theorem 3.8 in \cite{vdv}. $\Box$

 \bigskip Given a metric space $(X, D)$ let ${\sf Bd}(X)$ (respectively, ${\sf Comp}(X)$),  be the family of non empty bounded (respectively, compact) subsets of $X$. 
  \\For $S\subset X$ and $\delta > 0$ let 
 $\mathcal{N}_\delta(S) = \{x\in X : D(x, S) < \delta\}$; for  $X = E$ and $D$ the distance associated to $\norm \cdot \norm_{\un}$ (respectively, $\norm \cdot \norm_{\hun}$) this is $U_\delta(S)$ (respectively, $V_\delta(S)$). 
 Recall that the Hausdorff metric on  ${\sf Bd}(X)$ associated to the metric $D$ on $X$ is given by  $\mathcal{H}_D(S_1, S_2) = \inf\{\delta > 0 : S_1\subset \mathcal{N}_\delta(S_2) \text{ and } 
 S_2\subset \mathcal{N}_\delta(S_1)\} = \max\big\{\sup_{x\in S_1}D(x, S_2), \sup_{x\in S_2}D(x, S_1) \big\}$. 
 For $X = E$ and $D$ the distance associated to $\norm \cdot \norm_{\un}$ (respectively, the distance ${\sf D}_{\hun}$) the corresponding Hausdorff metric will be written  $\mathcal{H}_{\un}$ (respectively,  $\mathcal{H}_{\hun}$). From the first part of Lemma \ref{unifconv}, $\mathcal{H}_{\un} \leq \mathcal{H}_{\hun} \leq 2\mathcal{H}_{\un}$.
 
 \begin{prop}\label{contH1} With respect to either of the equivalent metrics  $\mathcal{H}_{\hun}$ or $\mathcal{H}_{\un}$,  the convex hull operator $S\mapsto\lb S\rb$ is Lipschitz continuous  from 
 ${\sf Bd}(E)$ to itself.
 \end{prop}
 
 \proof First, from Corollary \ref{diam}, $\lb S \rb$ is bounded if $S$ is bounded. Choose 
 $\delta > \mathcal{H}_{\hun}(S_1, S_2)$ then $S_1\subset V_\delta(S_2)\subset U_\delta(S_2)$; from $S_2\subset \lb S_2\rb$ we have 
 $\lb S_1\rb\subset \lb\, U_\delta(\lb S_2\rb)\,\rb$ and finally, from Lemma \ref{unifconv}, 
 $\lb\, U_\delta(\lb S_2\rb)\,\rb = U_\delta(\lb S_2\rb)\subset V_{2\delta}(\lb S_2\rb)$. We have shown that 
 $\lb S_1\rb\subset V_{2\delta}(\lb S_2\rb)$ and similarly $\lb S_2\rb\subset V_{2\delta}(\lb S_1\rb)$ from which 
 $\mathcal{H}_{\hun}(\lb S_1\rb, \lb S_2\rb) \leq 2 \mathcal{H}_{\hun}(S_1, S_2)$. $\Box$
 
 \bigskip Given a non empty max-plus convex $C$ subset of $E$ let ${\sf MPCC}_{\un}(C)$  be the family of non empty   max-plus convex and compact subsets  of $C$; 
 
 \begin{coro}\label{arhyperm+} 
 If $C$ is a non empty compact max-plus convex subset of $E$ then,  ${\sf MPCC}_{\un}(C)$ is an absolute retract. \end{coro} 
 
 \proof Since $C$ is compact it is, with respect to either of the metrics associated to the norms $\norm\cdot\norm_{\un}$ or 
 $\norm\cdot\norm_{\hun}$ a Peano continuum\footnote{A compact, connected and locally connected metric space} therefore, by Wojdislawki's Theorem \cite{wodj},  ${\sf Comp}(C)$ is an absolute retract; by Proposition \ref{hullcomp} and Proposition \ref{contH1}, $K\mapsto\lb K\rb$ is a continuous retraction from ${\sf Comp}(C)$ to  
 $\sf {\sf MPCC}_{\un}(C)$ which  is consequently an absolute retract. $\Box$
 
 \begin{coro}\label{arcompm+}
 A non empty compact and max-plus convex subset is an absolute retract.
 \end{coro}
 
 \proof Let $C$ be a non empty compact and maxplus convex subset of $E$. The map $x\mapsto\{x\}$ is continuous from $C$ to ${\sf CCMP}_{\un}(C)$. An arbitrary non empty max-plus  $K$ has a unique maximal element $\bigvee K$ and, for all $x\in C$, $\bigvee\{ x\} = x$; the map 
 $K\mapsto \bigvee K$ is continous and onto from ${\sf MPCC}_{\un}(C)$ to $C$. In conclusion, $C$ is a retract of   ${\sf MPCC}_{\un}(C)$. $\Box$
 
 \bigskip We will see in the next section that, in Corollary \ref{arcompm+}, the compactness assumption is superfluous, Theorem \ref{main}, but the proof of that result is somewhat more involved.  
 
 \bigskip If $C$ is a compact max-plus convex set $\bbR^n$ ($\un = (1, 1\cdots, 1)$),  it has been shown by  L. Bazylevych, M. Zarichnyi \cite{bazyzar} that ${\sf MPP}_{m}(C)$, the space of (non enmpty) max-plus polytopes of the form $\lb S\rb$ with $S\subset C$ and of cardinality at most $m$,  is an absolute retract. In that same paper they also show that, for a compact and metrizable space $\Omega$,  whose unit is taken to be the constant map $\omega\mapsto 1$, ${\sf MPCC}\big(\mathcal{C}(\Omega)\big)$ is an absolute retract homeomorphic to $l_2$.
 
 \medskip The compact topological space $\Omega$ of the  BBK Representation Theorem is unique up to homeomorphism; and it can be realized as a closed subspace of the unit ball of the dual space endowed with the weak topology (hence its compactness, by Alaoglu's Theorem); also, if $\big(E, \norm \cdot \norm_{\un}\big)$ is separable then the unit ball of the dual space is metrizable in the weak topology. The   BBK Representation Theorem combined with the  theorem of     L. Bazylevych and  M. Zarichnyi from \cite{bazyzar}   cited above  yields the following proposition. 
 
 \begin{prop}\label{arcmpp}
 Let $E$ be a Riesz space with unit $\un$. If $\big(E, \un\big)$ is complete and separable then ${\sf MPCC}_{\un}(E)$, the hyperspace of max-plus non empty compact subsets of $E$, is an absolute retract.
 \end{prop}
 
\section{Fixed points and selections}\label{topmax+conv}

\bigskip From Proposition \ref{quasiconmet} one has a max-plus version of Fan's Best Approximation Theorem, (A) of Proposition \ref{fanbest} below,   from which one has the max-plus version of Brouwer's Fixed Point Theorem. For the proof, in the context of geodesic spaces,   the reader is referred to \cite{midpt}, and to \cite{dugra} page 146 for the proof of the original Fan's Theorem in normed spaces. 

\begin{prop}\label{fanbest} Let $f : C\to E$ be a continuous function, with respect to ${\sf D}_{\hun}$, defined on a compact max-plus convex subset of $E$. If the metric  ${\sf D}_{\hun}$ is complete then the following hold: \\
{\bf (A)} {\bf (Fan's Best Approximation Theorem in max-plus)} There exists $x_0\in C$ such that, for all $y\in C$, 
${\sf D}_{\hun}(x_0, f(x_0)) \leq {\sf D}_{\hun}(y, f(x_0))$.\\
{\bf (B)} {\bf (Fan's Fixed Point Criteria in max-plus)} For $f$ to have a fixed point it is sufficient that, for all 
$x\in C$ such that $x\neq f(x)$, $C\cap \lb x, f(x) \rb$ contains a point other than $x$.
\end{prop}

One could relax the completness assumption on the metric ${\sf D}_{\hun}$ by assuming that $C$ is contained in a complete max-plus convex subset $X$ of $E$ and that $f(C)\subset X$. If $f$ takes its values in $C$ itself then, either (A) or (B) of Proposition \ref{fanbest},  implies that $f$ has a fixed point. This ``Brouwer's Fixed Point Theorem in max-plus'' follows also from the more general  ``Kakutani's Fixed Point Theorem in max-plus'', Theorem \ref{fanhi} below. \\
 Recall that the metric ${\sf D}_{\hun}$ is equivalent to the metric associated to the norm $\norm \cdot \norm_{\un}$ and  if there is on $E$ a complete Riesz norm then, for all $\un\in E_+$, $E_{\un}$ equipped with the norm $\norm \cdot \norm_{\un}$ is complete, \cite{frem} page 65.

\begin{thm}[Michael's Selection Theorem]\label{michael} If $X\subset E$ is a max-plus convex subset of $E$ which is complete with respect to ${\sf D}_{\hun}$ then, 
all lower semicontinuous maps $\Gamma : Y\to X$ with non empty closed
max-plus convex values defined a  paracompact space $Y$ have a continuous selection. 

\medskip Furthermore, if $A\subset Y$ is a closed set then any continuous selection of the resriction of $\Gamma$  to $A$ extends to a 
continuous selection of $\Gamma$.
\end{thm} 

\proof From  Theorem 3.4 in \cite{moitopconv}. $\Box$

\begin{thm}[Approximate selections for usc maps]\label{appscs} 
Let $X$ be a non empty max-plus convex subset of $E$, $Y$ a paracompact  topological space and $\Gamma : Y\to X$ an upper
semicontinuous map with non empty max-plus convex values. Then, for all $\delta > 0$, 
there exists a continuous map $f: Y\to X$ such that, for all $y\in Y$, $f(y)\in V_\delta(\Gamma y)$. 

\medskip Furthermore, if the values of $\Gamma$ are max-plus convex and compact then any neighborhood 
$\Theta\subset Y\times X$ of the graph of $\Gamma$ contains the graph of a continuous 
map $f: Y\to X$.
\end{thm}

\proof From  Theorem 3.5 in \cite{moitopconv}. $\Box$

\begin{thm}[Dungundji's Extension Theorem]\label{theodugu} 
Let $A\subset X$ be non empty closed subset of an arbitrary metric space $(X, d)$ and $f : A\to C$ a continuous map from $A$ to   an arbitrary non empty max-plus convex subset of $E$. Then, there exists a continuous maps  
$\hat{f} : X\to \lb f(A)\rb$ whose restriction to $A$ is $f$.
\end{thm}

\proof From Theorem 4.1 in \cite{moitopconv}. $\Box$
 
\bigskip Theorem \ref{theodugu} says that max-plus convex subsets (with respect to a given unit $\un$) of a Riesz space $E$ are, with respect to the topology induced by either of the norms $\norm \cdot \norm_{\un}$ or $\norm \cdot \norm_{\hun}$, absolute extensors for the class of metric spaces. 

\begin{thm}
\label{main} 
An arbitrary non empty max-plus convex subset of $E$, equipped with the metric topology associated to ${\sf D}_{\hun}$, is an
 absolute retract. 
\end{thm}

\proof A metrizable absolute extensor for the class of metric spaces is an absolute retract. $\Box$

\bigskip A map (single valued or multivalued)  $\Gamma : X \to Y$, where $Y$ is a toplological space, is a {\bf compact map} if there is a compact set $K\subset Y$ such that $\Gamma(Y)\subset K$. A set $X$ has {\bf the fixed point property} for a given class 
$\mathcal{M}$ of maps $\Gamma : X\to X$ if, for all $\Gamma\in\mathcal{M}$ there exists $x\in X$ such that $x\in\Gamma x$. 

\begin{thm}[Kakutani-Fan - Himmelberg's Theorem]\label{fanhi} An arbitrary non empty max-plus convex subset of $E$
has the fixed point property for upper semicontinuous  compact maps with closed  max-plus convex non empty values.
\end{thm}

\proof From Proposition \ref{toppolymax}, Lemma \ref{kakuprop} and Theorem 5.2 in \cite{moitopconv}. $\Box$

\bigskip Since a max-plus convex set is an absolute retract, Theorem \ref{fanhi} follows from the much harder Eilenberg-Montgomery Fixed Point Theorem, Corollary 7.4 in \cite{dugra} or from (iv) of Corollary 7.5 on the same page of which the following statement is a particular instance :  

\medskip\noindent{\it If $X$ is an absolute retract then, arbitrary compact upper semicontinous maps $S : X\to X$ with contactible values have a fixed point.} 

\bigskip From the  results of this section one can derive  tropical versions of some classical results from mathematical economics, existence of Nash equilibria, existence of equilibria for abstract economies, existence of maximal elements for a preference relation.  The question of their relevance in economics or  game theory is left open. 

\section{A few drawings in $\bbR^2$}
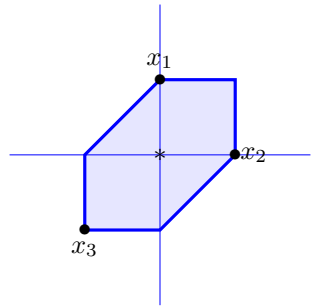
\begin{equation*}% [[x_1, x_2, x_3]] = D-ball
\begin{tikzpicture}[scale = 0.5]
\draw[color = white,  fill=blue!10 ] (0, 2) -- (2, 2) -- (2, 0) -- (0, -2) -- (-2, -2) -- (-2, 0) -- (0, 2) ;

\draw[color = blue, very thick]  (0, 2) -- (2, 2) -- (2, 0) -- (0, -2) -- (-2, -2) -- (-2, 0) -- (0, 2) ;
\draw[color = blue, very thin] (-4, 0)--(4, 0) ; 
\draw[color = blue, very thin] (0, 4)--(0, -4) ; 
\path[color = black] (0, 0) node(n) {$\ast$} ;
\path[color = black] (0, 2) node(n) {$\bullet$};
\path[color = black] (0, 2.5) node(n) {$x_1$};
\path[color = black] (2, 0) node(n) {$\bullet$};
\path[color = black] (2.5, 0) node(n) {$x_2$};
\path[color = black] (-2, -2) node(n) {$\bullet$};
\path[color = black] (-2, -2.5) node(n) {$x_3$};
\end{tikzpicture}
\end{equation*}

The closed unit  ${\sf D}_{\hun}$-ball in $\bbR^2$ about $\ast$; it is $\lb x_1, x_2, x_3\rb$.

\bigskip The next example shows that the metric ${\sf D}_{\hun}$ does not have the unique extension property. 

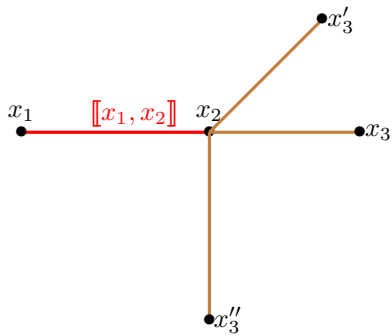
\begin{equation*}
\begin{tikzpicture}[scale = 0.5]
\draw[color = red, very thick] (-9, 3) -- (-4, 3) ; 
\path[color = black] (-9, 3) node(w) {$\bullet$};
\path[color = black] (-4, 3) node(w) {$\bullet$};
\path[color = black] (-9, 3.5) node(w) {$x_1$};
\path[color = black] (-4, 3.5) node(w) {$x_2$};
\path[color = red] (-6, 3.5) node(w) {$\lb x_1, x_2\rb$};
\draw[color = brown, very thick] (-4, 3) -- (0, 3) ;
\draw[color = brown, very thick] (-4, 3) -- (-1, 6) ;
\draw[color = brown, very thick] (-4, 3) -- (-4, -2) ;

\path[color = black] (0, 3) node(w) {$\bullet$};
\path[color = black] (0.5, 3) node(w) {$x_3$};
\path[color = black] (-1, 6) node(n) {$\bullet$};
\path[color = black] (-0.5, 6) node(w) {$x_3^\prime$};
\path[color = black] (-4, -2) node(w) {$\bullet$};
\path[color = black] (-3.5, -2) node(w) {$x_3^{\prime\prime}$};

\end{tikzpicture}
\end{equation*}

Three possible geodesic extensions of the geodesic $\lb x_1, x_2\rb$

\vspace{0.3cm}
$$\lb x_1, x_2\rb\subset \lb x_1, x_3\rb\quad\lb x_1, x_2\rb\subset \lb x_1, x^\prime_3\rb \text{ and } \lb x_1, x_2\rb\subset \lb x_1, x^{\prime\prime}_3\rb$$
 
  \medskip\centerline{$\lb x_1, x_3\rb$, $\lb x_1, x_3^\prime\rb$, $\lb x_1, x_3^{\prime\prime}\rb$ are geodesics.}

\subsection{A few max-plus polytopes in $\bbR^2$ and  a max-plus convex set that is not a polytope} 

In each case the max-plus polytope in question is the max-plus convex hull of the points labeled $x_1, x_2, x_2, \ldots$

\begin{equation*}% [[x_1, x_2, x_3]]
\begin{tikzpicture}[scale = 0.5]
\draw[color = blue, very thick] (-8, 0) -- (-5, 0) ; 
\draw[color = blue, very thick] (-5, 0) -- (-2, 3) ; 
\draw[color = blue, very thick] (-5, -4) -- (-5, 0) ;
\path[color = black] (-8, 0) node(w) {$\bullet$};
\path[color = black] (-2, 3) node(w) {$\bullet$};
\path[color = black] (-5, -4) node(w) {$\bullet$};
\path[color = black] (-8.5, 0) node(w) {$x_1$};
\path[color = black] (-2.2, 3.5) node(w) {$x_2$};
\path[color = black] (-4.5, -4) node(w) {$x_3$};
\end{tikzpicture}
\end{equation*}
This example shows that in $\bbR^2$ a max-plus convex polytopes with three extreme points (none of the three points is in the max-plus convex hull of the other two) can have empty interior and have topological dimension equal to $1$. 

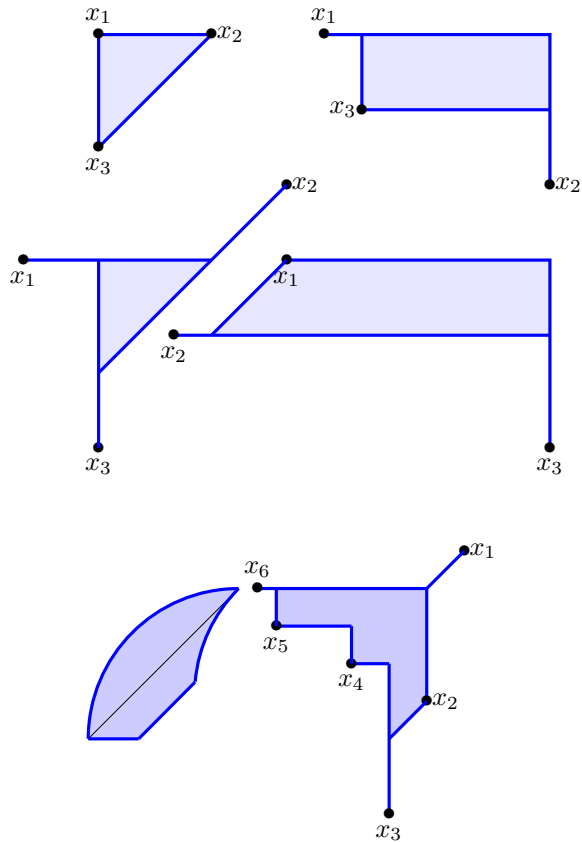
\begin{equation*}% [[x_1, x_2, x_3]]
\begin{tikzpicture}[scale = 0.5]
\draw[color = white,  fill=blue!10 ] (-8, 0) -- (-5, 0) -- (-8, -3) -- (-8, 0);
\draw[color = blue, very thick] (-8, 0) -- (-5, 0) -- (-8, -3) -- (-8, 0);
\path[color = black] (-8, 0) node(w) {$\bullet$};
\path[color = black] (-5, 0) node(w) {$\bullet$};
\path[color = black] (-8, -3) node(w) {$\bullet$};
\path[color = black] (-8, 0.5) node(n) {$x_1$};
\path[color = black] (-4.5, 0) node(n) {$x_2$};
\path[color = black] (-8, -3.5) node(s) {$x_3$};
\draw[color = white, fill=blue!10] (-1, 0) -- (-1, -2) -- (4, -2) -- (4, 0) -- (-1, 0) ;
\draw[color = blue, very thick] (-2, 0) -- (4, 0) -- (4, -4) ;
\path[color = black] (-2, 0) node(w) {$\bullet$};
\path[color = black] (-2, 0.5) node(w) {$x_1$};
\path[color = black] (4, -4) node(w) {$\bullet$};
\path[color = black] (4.5, -4) node(w) {$x_2$};
\draw[color = blue, very thick] (-1, 0) -- (-1, -2) -- (4, -2) ;
\path[color = black] (-1, -2) node(w) {$\bullet$};
\path[color = black] (-1.5, -2) node(w) {$x_3$};
\draw[color = white,  fill=blue!10 ] (-8, -6) -- (-5, -6) -- (-8, -9) -- (-8, -6);
\draw[color = blue, very thick] (-10, -6) -- (-5, -6);
\path[color = black] (-10, -6) node(w) {$\bullet$};
\path[color = black] (-8, -11) node(w) {$\bullet$};
\path[color = black] (-3, -4) node(w) {$\bullet$};
\path[color = black] (-10, -6.5) node(n) {$x_1$};
\path[color = black] (-2.5, -4) node(n) {$x_2$};
\path[color = black] (-8, -11.5) node(s) {$x_3$};
\draw[color = blue, very thick] (-8, -6) -- (-8, -11);
\draw[color = blue, very thick] (-3, -4) -- (-8, -9);
\draw[color = white,  fill=blue!10 ] (-3, -6) -- (4, -6)-- (4, -8) -- (-5, -8) -- (-3, -6);
\draw[color = blue, very thick]  (-3, -6) -- (4, -6) -- (4, -11);
\path[color = black] (4, -11) node(w) {$\bullet$};
\path[color = black] (4, -11.5) node(s) {$x_3$};
\path[color = black] (-3, -6) node(w) {$\bullet$};
\path[color = black] (-3, -6.5) node(w) {$x_1$};
\draw[color = blue, very thick]  (-3, -6) -- (-5, -8) ;
\draw[color = blue, very thick]  (-6, -8) -- (4, -8) ;
\path[color = black] (-6, -8) node(w) {$\bullet$};
\path[color = black] (-6, -8.5) node(w) {$x_2$};
\end{tikzpicture}
\end{equation*}

\begin{equation*}
\begin{tikzpicture}[scale = 0.5]

\draw[color = blue,   fill = blue!20, very thick] (0, 0) arc (90:180:4cm);
\draw[color = white, fill = blue!20] (-4, -4) -- (-2.66 , -4) --  (-1.17, -2.51) -- (0, 0)-- cycle ;
\draw[color = black] (0, 0) -- (-4, -4);
\draw[color = blue, fill = white, , very thick] (0, 0) arc (135:175:4cm);
\draw[color = blue, very thick]  (-4, -4) -- (-2.66 , -4) ;
\draw[color = blue, very thick]   (-1.159, -2.49) -- (-2.66 , -4) ;
\draw[color = white, fill = blue!20] (6, 1) -- (5, 0) -- (5, -3) -- (4, -4) -- (4, -2) -- (3, -2) -- (3, -1) -- (1, -1) -- (1, 0)--(0.5, 0) -- (5, 0) -- (6, 1) ;
\path[color = black] (6, 1) node(w) {$\bullet$};
\path[color = black] (6.5, 1) node(w) {$x_1$};
\draw[color = blue, very thick]  (6, 1) -- (5, 0) ; 
\draw[color = blue, very thick]  (5, 0) -- (5, -3) ;
\path[color = black] (5, -3) node(w) {$\bullet$};
\path[color = black] (5.5, -3) node(w) {$x_2$};
\draw[color = blue, very thick]  (5, -3) -- (4, -4) ;
\draw[color = blue, very thick] (4, -6) -- (4, -2) ;
\path[color = black] (4, -6) node(w) {$\bullet$};
\path[color = black] (4, -6.5) node(w) {$x_3$};
\draw[color = blue, very thick] (4, -2) - - (3, -2) ;
\path[color = black] (3, -2) node(w) {$\bullet$};
\path[color = black] (3, -2.5) node(w) {$x_4$};
\draw[color = blue, very thick] (3, -2) -- (3, -1); 
\draw[color = blue, very thick] (3, -1) -- (1, -1);
\path[color = black] (1, -1) node(w) {$\bullet$};
\path[color = black] (1, -1.5) node(w) {$x_5$};
\draw[color = blue, very thick] (1, 0)--(1, -1);
\draw[color = blue, very thick] (0.5, 0) -- (5, 0);
\path[color = black] (0.5, 0) node(w) {$\bullet$};
\path[color = black] (0.5, 0.5) node(n) {$x_6$};
\end{tikzpicture}
\end{equation*}

In $\bbR^n$ a max-plus segment is  piecewise linear, it is made of at most $n$ affine segments. A max-plus polytope is a contractible finite union of affine convex polytopes, and consequently an absolute retract; it is also a contractible simplicial  complex. What is the structure of max-plus polytopes in infinite dimensional Riesz spaces?

\end{document}